\documentclass[11pt]{article}
\usepackage{amsmath}
\usepackage{amsthm}
\usepackage{amsfonts}
\usepackage{amssymb}
\usepackage{mathtools}
\usepackage[numbers,sort&compress]{natbib}

\usepackage[draft]{todonotes}   

\usepackage{cases}
\usepackage{xcolor}

\newcommand{\dt}{\partial_t}
\newcommand{\dv}{\mathrm{div}\,}

\newcommand{\idx}[0]{\,dx}

\newcommand{\abs}[2]{\bigl| #1 \bigr|^{#2}}
\newcommand{\norm}[2]{\Arrowvert #1 \Arrowvert_{#2}}

\newcommand{\Lnorm}[1]{L^{#1}(\Omega)}
\newcommand{\Hnorm}[1]{H^{#1}(\Omega)}

\renewcommand{\vec}[1]{\mathbf{#1}}

\newcommand{\mrm}[1]{\mathrm{#1}}

\newtheorem{thm}{Theorem}[section]

\theoremstyle{remark}

\numberwithin{equation}{section}
\allowdisplaybreaks
\title{Well-posedness of Hibler's dynamical sea-ice model}

\author{Xin Liu\footnote{Weierstrass-Institut
f\"ur Angewandte Analysis und Stochastik,
Leibniz-Institut im Forschungsverbund Berlin, Berlin Germany. Isaac Newton Institute for Mathematical Sciences, University of Cambridge, Cambridge CB3 0EH, UK.
Email: stleonliu@gmail.com\, and \, stleonliu@live.com}
\qquad Marita Thomas\footnote{Weierstrass-Institut
f\"ur Angewandte Analysis und Stochastik,
Leibniz-Institut im Forschungsverbund Berlin, Berlin Germany. Email: thomas@wias-berlin.de} 
\qquad Edriss S. Titi\footnote{Department of Applied Mathematics and Theoretical Physics, University of Cambridge, Cambridge CB3 0WA, UK. Department of Mathematics, Texas A{\&}M University, College Station,  TX 77840, USA.
		Department of Computer Science and Applied Mathematics, Weizmann Institute of Science, Rehovot 76100, Israel. Email: {Edriss.Titi@damtp.cam.ac.uk}\, and \, {titi@math.tamu.edu}}}

\begin{document}
\allowdisplaybreaks
\maketitle

\begin{abstract}
	This paper establishes the local-in-time well-posedness of solutions to an approximating system constructed by mildly regularizing the dynamical sea-ice model of {\it W.D. Hibler, Journal of Physical Oceanography, 1979}. Our choice of regularization has been carefully designed, prompted by physical considerations, to retain the original coupled hyperbolic-parabolic character of Hibler's model. Various regularized versions of this model have been used widely for the numerical simulation of  the circulation and thickness of the Arctic ice cover. However, due to the singularity in the ice rheology, the notion of solutions to the original model is unclear. Instead, an approximating system, which captures current numerical study, is proposed. The well-posedness theory of such a system provides a first-step groundwork in both numerical study and future analytical study. 
\end{abstract}

\tableofcontents

\section{Introduction}

\subsection{The sea-ice dynamic-thermodynamic model}
Global climate changes, especially global warming, have large impact on the Arctic sea-ice, which has, in return, determining effects on not only global climate but also the local and global ecosystem, human activities etc. (see e.g., \cite{seaice-second-edition}). 
If the problem is statically determinate, as pointed out in \cite{Schreyer2001}, a sea-ice dynamical model based on the viscous-plastic rheology was introduced in \cite{Hibler1979}, where the thickness of ice plays an essential role in the thermodynamics, and characterizes the strength of the ice interaction (i.e., ice rheology). The velocity of sea-ice $ \vec u $ is described by two-dimensional momentum balance equations, where the viscosity effect is characterized by a viscous-plastic rheology, and the strength of viscosity depends on the thickness of ice. The mean ice thickness $ h $ and the compactness of ice $ A $ are described by two continuity equations with thermodynamic source terms. That is, with a simplified ice rheology (see \eqref{def:ice-rheology-app}, below), the above quantities are governed by the coupled system,
\begin{subequations}\label{eq:sea-ice-origin}
\begin{gather}
	m (\dt \vec u + \vec u \cdot \nabla \vec{u}) + \nabla p = \dv \mathbb S + \mathcal F, \\
	\dt h + \dv (h\vec u) = \mathcal S_h, \\
	\dt A + \dv (A\vec u) = \mathcal S_A + A \dv \vec u \cdot \chi_{\lbrace A \geq 1 \rbrace} ,
\end{gather}	
\end{subequations}
with
\begin{subequations}
\begin{align}
	\text{ice mass} ~ m & := \rho_\mrm{ice} h, \label{def:ice-mass} \\
	\text{pressure} ~  p & := c_p h \exp (c_a A), \label{def:pressure} \\
	\text{viscoplastic stress} ~  \mathbb S & := p \dfrac{\nabla u + \nabla u^\top}{\abs{\nabla u + \nabla u^\top}{}}	 + p \dfrac{\dv \vec u \mathbb I_2}{\abs{\dv \vec u}{}}, \label{def:ice-rheology} \\ 
	\mathcal F & := - m \eta \vec u^{\perp} + \vec \tau_\mrm a + \vec \tau_\mrm w, \label{def:forcing}\\
	\text{air flow stress} ~ \vec \tau_\mrm a & := \rho_\mrm a C_\mrm a \abs{\vec U_\mrm g}{} ( \vec U_\mrm g \cos \phi + \vec U_\mrm g^\perp \sin \phi), \label{def:air-stress} \\
	\text{water flow stress} ~ \vec \tau_\mrm w & := \rho_\mrm w C_\mrm w \abs{\vec U_\mrm w - \vec u}{} \lbrack ( \vec U_\mrm w - \vec u) \cos \theta {\nonumber}  \\
	 & \qquad\qquad + ( \vec U_\mrm w - \vec u)^\perp \sin \theta \rbrack, \label{def:water-stress} \\
	 \mathcal S_h & := \bigl\lbrack f({h}/{A}) A + (1-A) f(0)\bigr\rbrack \cdot \chi_{\lbrace h > 0 \rbrace}, \label{def:h-source} \\
	\mathcal S_A & :=  \bigl({\bigl(f(0)\bigr)^+}/{h_0}\bigr)(1-A) + ( - A/(2h) ) \cdot ( \mathcal S_h )^-  . \label{def:A-source}
\end{align}
Here $ \chi_{\lbrace h > 0 \rbrace }, \chi_{\lbrace A \geq 1 \rbrace }$ are the characteristic functions of sets $ \lbrace h > 0 \rbrace ,\lbrace A \geq 1 \rbrace  $, defined by
\begin{equation}\label{def:cha-fun}
\chi_{\lbrace h > 0 \rbrace} =
\begin{cases}
	1 & h > 0, \\
	0 & h \leq 0,
\end{cases}
~~~~
 \chi_{\lbrace A \geq 1 \rbrace } = \begin{cases}
	1 & A \geq 1, \\
	0 & A < 1,
\end{cases}
\end{equation}
respectively.
\end{subequations}
In addition, $ \vec v^\perp = (-v_2, v_1)^\top $ for any vector $ \vec v = (v_1, v_2)^\top $; $ \rho_\mrm{ice}, \rho_\mrm{a}, \rho_\mrm{w} $ represent the density of ice, air, and water, respectively; $ c_p, c_a, C_\mrm a,  C_\mrm w $ are the thermodynamic constants; and $ \vec U_\mrm g, \vec U_\mrm w, \phi, \theta $ denote the velocity and stress angle of the air and the water, which, for simplicity of presentation, are assumed to be constant in this paper.

System \eqref{eq:sea-ice-origin} is used to study simulation the evolution of sea-ice in numerical study. For instance, the model successfully reproduces many of the observed features of the circulation and thickness of the Arctic ice cover in \cite{Hibler1979}. Hibler's model of sea-ice dynamics is the foundation for further model developments, including the elastic-viscous-plastic sea-ice dynamics model as in \cite{Hunke2001}, the Maxwell elasto-brittle rheology model as in \cite{DWSL2016}, and models with leads, ridging or tensile failure as in \cite{SSMCK2006,WF2012}. See \cite{Pritchard2001} for a summary of classical models with different descriptions of ice distribution, and  \cite{BFLM2013,Hunke1997,Mehlmann2017a,Parkinson1979,RBOM2016,Herman2016,TFW2013,CKEP1998} and the references therein for further model development and computational investigation. 


Despite the steady advances that have been made in 
 the modeling and simulation of sea ice, mathematical analysis is much less developed in this field. It is the main objective of this paper to provide rigorous mathematical analysis of Hibler's model as a first step in this direction and as a basis for further investigation of next-generation sea-ice rheologies.


In particular, the fundamental problem of well-posedness of solutions to system \eqref{eq:sea-ice-origin} is widely open, which is related to the validity of the model as pointed out by \cite{Schreyer2001}. In \cite{GK1995,Gray1999}, the authors study the loss of hyperbolicity of the linearized system of the original Hibler's model around divergent flows, and show that the system is ill-posed. This work is further discussed in \cite{Dukowicz1997,LHMJ2007,Guba2013,ST2015} from various perspectives. 
These studies do not contradict the local well-posedness result established in this paper, for the following reasons:
\begin{itemize}
	\item Instead of the original Hibler's model, we consider a regularization of Hibler's model, which, as shown below in the paper, preserves the parabolicity of the momentum equations in a Sobolev space with enough regularity. This is very different from the hyperbolic equations considered in the ill-posedness studies;
	\item Instead of linear analysis, we consider the nonlinear well-posedness theory for regularized Hibler's model in the Hadamard sense, including the existence, the uniqueness, and the continuous dependency on initial data in the suitable Sobolev space as shown in Theorem \ref{thm:well-posedness-app-2}. 
\end{itemize}


We would like to point out that the main challenge in establishing the well-posedness theory is the singularity arising in the stress tensor \eqref{def:ice-rheology} when $ |\nabla \vec u| \rightarrow 0^+ $. In fact, among the numerical investigations, such singularity is usually truncated, i.e. regularized, by replacing it with its strictly positive approximation (e.g., $ \max\lbrace|\nabla \vec u|, \Delta_\mathrm{min} \rbrace $ or $ \sqrt{|\nabla \vec u|^2 + \varepsilon^2} $).


Notably, we would like to point out an investigation of very singular diffusion equations in \cite{Giga1998,Giga2010}, where the authors discuss the notion of solutions to
$$
\dt u = \dv \biggl( \dfrac{\nabla u}{|\nabla u|}\biggr).
$$
Similarly, the positive one-homogeneity of the potential related to \eqref{def:ice-rheology} calls for a subdifferential formulation of the problem, however set in the Eulerian frame. We leave such investigation to our future study.

In this paper, due to the obstacles mentioned above, we propose to study the following regularized approximating problem of \eqref{eq:sea-ice-origin}: for $ \varepsilon, \omega \in (0,1) $,
\begin{subequations}\label{eq:sea-ice-app-2}
    \begin{gather}
        m(\dt \vec u + \vec u \cdot \nabla \vec u) + \nabla p = \dv \mathbb S_{\varepsilon} + \mathcal F, \\
        \dt h + \dv(h \vec u) = \mathcal S_{h,\omega}, \\
        \dt A + \dv(A \vec u) = \mathcal S_{A,\omega} + A \dv \vec u \cdot\chi^\omega_A,
    \end{gather}
\end{subequations}
where $ m $, $ p $, and $ \mathcal F$ are as in \eqref{def:ice-mass}, \eqref{def:pressure},  and \eqref{def:forcing}, respectively, and
\begin{subequations}\label{eq:sea-ice-app-2-rh}
\begin{align}
	\mathbb S_\varepsilon = \mathbb S_\varepsilon(p,\nabla \vec u) := & p \dfrac{\nabla \vec u + \nabla \vec u^\top}{\sqrt{\abs{\nabla \vec u + \nabla \vec u^\top}{2}+\varepsilon^2}}+ p \dfrac{\dv \vec u \mathbb I_2}{\sqrt{\abs{\dv \vec u}{2}+\varepsilon^2}}, \label{def:ice-rheology-app-1} \\
    \mathcal S_{h,\omega} := & \bigl\lbrack f(h/(A+\omega)) A + (1-A) f(0) \bigr\rbrack \chi_{\lbrace h>0 \rbrace}, \label{def:h-source-app-2} \\
    \mathcal S_{A,\omega} := & \dfrac{\bigl(f(0)\bigr)^+}{h_0} (1-A) -\dfrac{A}{2h}\cdot \dfrac{\sqrt{|\mathcal S_{h,\omega}|^2 + \omega^2} - \mathcal S_{h,\omega}}{2}, \label{def:A-source-app-2}
    \\
	\chi^\omega_A  := & 1 -  \dfrac{(1-A)^+}{(1-A)^+ + \omega}.\label{def:cha-A-app}
\end{align}
\end{subequations}

 To be more precise, we will establish the local in time well-posedness of strong solutions to \eqref{eq:sea-ice-app-2} in domain $ \Omega := \mathbb T^2 \subset \mathbb R^2 $:
 \begin{thm}\label{thm:well-posedness-app-2}
 	Consider initial data \begin{equation}\label{cnts:initial-1} (\vec u, h, A)|_{t=0} = (\vec u_\mrm{in}, h_\mrm{in}, A_\mrm{in}) \in (\Hnorm{3})^3 \end{equation} to system \eqref{eq:sea-ice-app-2}, satisfying
 	\begin{equation}\label{cnts:initial-2}
 	0 < \underline h \leq h_\mrm{in} \leq \overline h <\infty, \quad \text{and} \quad 0 \leq A_\mrm{in} \leq 1.
 	\end{equation}
 	In addition, we assume that
 	\begin{equation}\label{cnts:f}
	\begin{gathered}
		\underline f \leq f \leq \overline f, \\
		\abs{f'}{} + \abs{f''}{} + \abs{f'''}{} \leq M_f,
	\end{gathered}
	\end{equation}
	for some constants $ \underline f, \overline f \in \mathbb R $, $ M_f \in (0,\infty) $.
 	Then there exists a unique strong solution $ (\vec u, h, A) $ to system \eqref{eq:sea-ice-app-2} in $ [0,T]\times \Omega $, for some $ T \in (0,\infty) $ depending on initial data, with
\begin{equation}\label{limit:regularity}
    \begin{gathered}
    \vec u \in L^\infty(0,T;H^3(\Omega)) \cap L^2(0,T;H^4(\Omega)), \\
    h, A \in L^\infty(0,T;H^3(\Omega)), \\
    \dt\vec u, \dt h, \dt A \in L^\infty(0,T;L^2(\Omega)),
\end{gathered}
\end{equation}
and
\begin{equation}\label{limit:est-regularity}
\begin{gathered}
\norm{\vec  u, h, A}{L^\infty(0,T;H^3(\Omega))} + \norm{\vec u}{L^2(0,T;H^4(\Omega))} \\
    + \norm{\dt \vec u,\dt h, \dt A}{L^\infty(0,T;L^2(\Omega))} \leq \mathfrak C_\mrm{in}, \\
     0 \leq A \leq 1, \qquad 0 < \dfrac{1}{4} \underline h \leq h \leq 4 \overline h,
    \end{gathered}
\end{equation}
 where $ \mathfrak C_\mrm{in} \in (0,\infty) $ is some positive constant depending only on initial data. Moreover, the solution is stable with respect to perturbation of initial data.
 \end{thm}

Now, let us explain our strategy. Instead of directly constructing solutions to system \eqref{eq:sea-ice-app-2}, we consider another regularized system, parametrized by $ (\mu,\lambda,\iota,\nu) \in (0,1)^4 $:
%
%
%
%
\begin{subequations}\label{eq:sea-ice-app}
\begin{gather}
	m (\dt \vec u + \vec u \cdot \nabla \vec{u}) + \nabla p = \dv \mathbb S_{\varepsilon, \mu, \lambda} -  \iota \Delta^2 \vec{u} + \mathcal F, \label{eq:sea-ice-app=1} \\
	\dt h + \dv (h\vec u) = \mathcal S_{h,\omega,\nu}, \label{eq:sea-ice-app=2}\\
	\dt A + \dv (A\vec u) = \mathcal S_{A,\omega,\nu} + A \dv \vec u \cdot \chi^\omega_A, \label{eq:sea-ice-app=3}
\end{gather}
\end{subequations}
where $ m $, $ p $, $ \mathcal F $, and $ \chi_A^\omega $ are as in \eqref{eq:sea-ice-app-2} and \eqref{eq:sea-ice-app-2-rh}, and
\begin{subequations}
\begin{align}
	\mathbb S_{\varepsilon, \mu, \lambda} & := \mathbb S_\varepsilon + \mathbb S_{\mu,\lambda}, \label{def:ice-rheology-app}  \\
	\mathcal S_{h,\omega,\nu} & := \bigl\lbrack f({h^+}/{(A^+ + \omega )}) A + (1-A) f(0)\bigr\rbrack \cdot \chi^\nu_h, \label{def:h-source-app} \\
	\mathcal S_{A,\omega,\nu} & :=  \dfrac{\bigl(f(0)\bigr)^+}{ h_0+\nu}(1-A)  - \dfrac{A}{2h^+ + \nu}  \cdot \dfrac{\sqrt{|\mathcal S_{h,\omega,\nu}|^2 + \omega^2} - \mathcal S_{h,\omega,\nu}}{2}  , \label{def:A-source-app}
\end{align}
with $ \mathbb S_\varepsilon $ as in \eqref{eq:sea-ice-app-2} and
\begin{align}
	 \mathbb S_{\mu,\lambda} & := \mu (\nabla \vec u + \nabla \vec u^\top) + \lambda \dv \vec u \mathbb I_2, \label{def:fluid-rheology} \\
	\chi^\nu_h & := \dfrac{h^+}{h^+ + \nu},  \label{def:cha-h-app}.
\end{align}
\end{subequations}
We will construct solutions to system \eqref{eq:sea-ice-app} through a contraction mapping argument. That is, we consider a ``linearization" of \eqref{eq:sea-ice-app}, and establish a contraction mapping with respect to $ L^2 $ topology with bounds in a smooth function space. Then with a uniform-in-$(\mu,\lambda,\iota,\nu)$ estimate, we will be able to pass the limit $ (\mu,\lambda,\iota,\nu) \rightarrow (0^+,0^+,0^+,0^+)$, and eventually construct the strong solution to \eqref{eq:sea-ice-app-2}. The proof of Theorem \ref{thm:well-posedness-app-2} is then finished by showing the uniqueness and continuous dependency on the initial data. We would like to mention that the key ingredient in establishing the well-posedness of solutions involves showing the monotonicity of $ \mathbb S_\varepsilon(\cdot) $ in $ \nabla \vec u $, which is not trivially obvious
due to the fact that $ \mathbb S_\varepsilon (\cdot) $ is nonlinear in $ \nabla\vec u$. In particular, we will require the inequality of the type
$$
\bigl( \mathbb S_\varepsilon(p_1, \nabla \vec u_1) - \mathbb S_\varepsilon(p_2, \nabla \vec u_2) \bigr) : \bigl( \nabla \vec u_1 - \nabla \vec u_2 \bigr) \gtrsim | \nabla ( \vec u_1 - \vec u_2 ) |^2 + \cdots.
$$
We successfully establish this inequality by writing $  \mathbb S_\varepsilon(p_1, \nabla \vec u_1) - \mathbb S_\varepsilon(p_2, \nabla \vec u_2) $ in a symmetric form (see \eqref{id:monotone}, below).

We would like to make some remarks before going into details of the proof. Our ice rheology \eqref{def:ice-rheology-app-1} is a simplified version of the one from \cite{Hibler1979}. For some technical reasons, we are not sure whether Theorem \ref{thm:well-posedness-app-2} will apply to the original ice rheology from \cite{Hibler1979}.
We have not successfully established a proper uniform-in-$\varepsilon $ estimates of the solutions to \eqref{eq:sea-ice-app-2}. Therefore, we have not yet been able to establish a proper notion of solutions to the original system \eqref{eq:sea-ice-origin}. However, our approximation \eqref{eq:sea-ice-app-2} agrees with the most common numerical approaches to \eqref{eq:sea-ice-origin}, which, as we explain before, is restricted to a truncated ice rheology. Thus, in this sense, our analytical results provide a solid ground for current numerical schemes of \eqref{eq:sea-ice-origin}. Another issue is that we only consider the case when $ h_\mrm{in} \geq \underline h > 0 $, i.e., there is no absence of ice in the domain of study. To carry out the limit $ \underline h \rightarrow 0^+ $, more comprehensive {\it a priori} estimates are required. We leave this to future study.

Recently we have learnt an independent study \cite{BDHH-2021} by Brandt, Disser, Haller-Dintelmann, and Hieber about similar model, 
where the domain boundary and boundary conditions are taken into account instead of periodic domain. 
It is worth pointing out that in our regularized system \eqref{eq:sea-ice-app-2} the governing equations for the evolution of $ h $ and $ A $ remains hyperbolic, and therefore  system  \eqref{eq:sea-ice-app-2} is a mixed type system, while the regularized system in \cite{BDHH-2021} is parabolic in all its components. 
 In particular, due to the hyperbolicity, system \eqref{eq:sea-ice-app-2} is expected to have a completely different long-time dynamics than those investigated in \cite{BDHH-2021}. Moreover, the additional dissipation introduced in \cite{BDHH-2021} allows the authors to show global existence for small initial data.

This paper is organized as follows. In the next subsection, we will summarize some notations used in this paper. In Section \ref{sec:app-scheme}, we will detail the approximation scheme to \eqref{eq:sea-ice-app}. In Section \ref{sec:first-layer-app}, we establish the well-posedness of solutions to \eqref{eq:sea-ice-app} via a contraction mapping argument. Finally in Section \ref{sec:well-posedness}, we establish the uniform-in-$(\mu,\lambda,\iota,\nu)$ estimates, and pass to the limit $ (\mu,\lambda,\iota,\nu) \rightarrow (0^+,0^+,0^+,0^+) $ to show the existence of solutions to \eqref{eq:sea-ice-app-2}. The well-posedness of solutions is then established in Section \ref{subsec:well-posedness-app-2}

\subsection{Notations}

We use $ L^p(\cdot) $ and $ H^s(\cdot) $ to denote the standard Lebesgue and Sobolev spaces, respectively. For any functional space $ \mathcal X $ and functions $ \psi, \phi, \cdots $, we denote by
$$
\norm{\psi,\phi,\cdots}{\mathcal X} := \norm{\psi}{\mathcal X} + \norm{\phi}{\mathcal X}  + \cdots.
$$
In addition,
$$
	\psi^+ := \begin{cases}
		\psi &\text{if} ~ \psi \geq 0,\\
		0 &\text{if} ~ \psi < 0,
	\end{cases} \qquad \psi^- = \psi^+ - \psi.
$$

Let $ \partial \in \lbrace \partial_x, \partial_y \rbrace $. For any multi-index $ (\alpha_1, \alpha_2)\in (\mathbb Z^+)^2 $, denote by $ \partial^\alpha := \partial_x^{\alpha_1} \partial_y^{\alpha_2} $ with $ \alpha = \alpha_1 + \alpha_2 $. Throughout this paper, we use the notation $ X \lesssim Y $ to represent $ X \leq C Y $ for some generic constant $ C \in (0,\infty) $, which may be different from line to line. We use $ C_{a,b,\cdots} $ to emphasize the dependency on the quantities $ a, b, \cdots $. In addition, by $ \mathcal H (\cdots) $, it represents a generic bounded function of the arguments.

\section{An approximation scheme to solve \eqref{eq:sea-ice-app}}\label{sec:app-scheme}

\subsection{A ``linearization'' of \eqref{eq:sea-ice-app}}

Given $ \vec u^o $, assumed to be smooth enough, we consider first the following coupled hyperbolic system
\begin{subequations}
\begin{gather}
	\dt h_\mrm{m} + \dv(h_\mrm{m} \vec u^o) = \mathcal S_{h_\mrm{m},\omega,\nu}, \label{eq:h-lin} \\
	\dt A_\mrm{m} + \dv(A_\mrm{m} \vec u^o) = \mathcal S_{A_\mrm{m},\omega, \nu} + A_\mrm{m} \dv \vec u^o \cdot \chi^\omega_{A_\mrm{m}}, \label{eq:A-lin}
\end{gather}
where $ \mathcal S_{h_\mrm{m},\omega, \nu}$, $ \mathcal S_{A_\mrm{m},\omega,\nu} $, $ \chi_{h_\mrm m}^\omega $, and $ \chi_{A_\mrm m}^\omega $ are defined as in \eqref{def:h-source-app}, \eqref{def:A-source-app}, \eqref{def:cha-h-app}, and \eqref{def:cha-A-app}, with $ h $ and $ A $ replaced by $ h_\mrm m $ and $ A_\mrm m $, respectively. Here we use the subscript $ \mrm m $ (short for `mapping') and the superscript $ o $ (short for `origin') to label outputs and inputs in our contraction mapping.

We claim that, at least locally in time, there exists a unique solution $ (h_\mrm m, A_\mrm m) $ to \eqref{eq:h-lin} and \eqref{eq:A-lin} with proper initial data, for smooth enough $ \vec u^o $. $ (h_\mrm m,A_\mrm m) $ can be arbitrarily regular, provided that $ \vec u^o $ and initial data are regular enough. We leave the investigation of the regularity of $ (h_\mrm m, A_\mrm m) $ in the subsequent sections.

We remark that such claims follow from the standard well-posedness theory of hyperbolic equations (see, e.g., \cite{Majda1984}). Hence the proof is omitted.

Let $ (h_\mrm m, A_\mrm m ) $ be the solution to \eqref{eq:h-lin} and \eqref{eq:A-lin} as above, and consider the following equation:
\begin{equation}\label{eq:u-lin}
		\rho_\mrm{ice} h_\mrm m \dt \vec u_\mrm m + \iota \Delta^2 \vec u_\mrm m = - \rho_\mrm{ice} h_\mrm m \vec u^o \cdot \nabla \vec u^o - \nabla p_\mrm m + \dv \mathbb S_{\varepsilon,\mu, \lambda,\mrm m} + \mathcal F_\mrm m,
\end{equation}
\end{subequations}
where $ p_\mrm{m} $,  $ S_{\varepsilon,\mu, \lambda,\mrm m} $, and $ \mathcal F_\mrm m $ are defined as in \eqref{def:pressure}, \eqref{def:ice-rheology-app}, and \eqref{def:forcing}, with $ h $, $ A $, and $ u $ replaced by $ h_\mrm m $, $ A_\mrm{m} $, and $ \vec u^o $, respectively.

To solve the linear equation \eqref{eq:u-lin} by, e.g., a Galerkin method, one will need to deal with the possible degeneracy of $ h_\mrm{m} $. For this, we subsequently show that for $ \vec u^o $ smooth enough, with appropriate initial data, $ h_\mrm{m} $ and $ A_\mrm{m} $ satisfy certain non-degeneracy property.

\subsection{Non-negativity and uniform bound of $ A_\mrm m $: $ 0 \leq A_\mrm m \leq 1 $}\label{subsec:pointwise-A-app}

In this subsection, we show that $ 0 \leq A_\mrm m \leq 1 $ for a smooth enough $ \vec u^o $. In fact, we only require that
\begin{equation}\label{rg:24-feb-3}
	\dv \vec u^o \in L^1(0,T;L^\infty(\Omega)),
\end{equation}
for some $ T > 0 $.

\vspace{5mm}

{\bf\noindent Non-negativity of $ A_\mrm m $:}

Taking the $ L^2 $-inner product of \eqref{eq:A-lin} with $ \bigl( - A_{\mrm m}^- \bigr) $ leads to, after applying integration by parts in the resultant
\begin{equation}\label{eq:positivity-A-app}
\begin{gathered}
	\dfrac{1}{2} \dfrac{d}{dt} \norm{A_\mrm m^-}{\Lnorm{2}}^2 = \int \biggl( \dfrac{1}{2} - \dfrac{(1-A_\mrm m)^+}{(1-A_\mrm m)^+ + \omega} \biggr) \dv \vec u^o \abs{A_\mrm{m}^-}{2} \idx \\
	+ \int \underbrace{ \mathcal S_{A_\mrm m ,\omega, \nu} ( - A_\mrm{m}^- )}_{\leq 0} \idx \lesssim \norm{\dv \vec{u}^o}{\Lnorm{\infty}}\norm{A_\mrm{m}^-}{\Lnorm{2}}^2.
\end{gathered}
\end{equation}
Therefore, applying Gr\"onwall's inequality to \eqref{eq:positivity-A-app} yields
\begin{equation*}
	\norm{A_\mrm m ^-}{\Lnorm{2}}^2 \leq e^{C\int_0^t \norm{\dv \vec u^o(s)}{\Lnorm{\infty}}\,ds} \norm{A_\mrm{in}^-}{\Lnorm{2}}^2 = 0,
\end{equation*}
which implies
\begin{equation*}
	A_\mrm{m} \geq 0.
\end{equation*}

\vspace{5mm}

{\bf\noindent Non-negativity of $ 1 - A_\mrm m $:}

Consider the following equation for $ 1 - A_\mrm m $, derived from \eqref{eq:A-lin}:
\begin{equation}\label{eq:1-A-lin}
	\dt (1-A_\mrm m ) = - \mathcal S_{A_\mrm m, \omega, \nu} - \vec u^o \cdot \nabla (1-A_\mrm m) + A_\mrm m \dv \vec u^o \dfrac{(1-A_\mrm m)^+}{(1-A_\mrm m)^+ + \omega}.
\end{equation}
As before, taking the $ L^2 $-inner product of \eqref{eq:1-A-lin} with $ \bigl\lbrack - (1-A_\mrm m)^- \bigr\rbrack $, after applying integration by parts, leads to
		\begin{gather*}
			 \dfrac{1}{2} \dfrac{d}{dt} \norm{(1-A_\mrm{m})^-}{\Lnorm{2}}^2 = \int \biggl( \dfrac{1}{2}   \dv \vec u^o \abs{(1-A_\mrm m)^-}{2}+ \underbrace{ \mathcal S_{A_\mrm m, \omega, \mu} (1-A_\mrm m)^- }_{\leq 0} \biggr) \idx \\
			 - \int \underbrace{A_\mrm m \dv \vec u^o \dfrac{(1-A_\mrm m)^+}{(1-A_\mrm m)^+ + \omega} (1-A_\mrm m)^-}_{= 0} \idx,
		\end{gather*}
which yields
\begin{equation}\label{eq:upper-A-app}
	\dfrac{d}{dt} \norm{(1-A_\mrm{m})^-}{\Lnorm{2}}^2
			 \lesssim \norm{\dv \vec u^o}{\Lnorm{\infty}} \norm{(1-A_\mrm{m})^-}{\Lnorm{2}}^2.
\end{equation}
Then as before, after applying Gr\"onwall's inequality to \eqref{eq:upper-A-app}, one can conclude
\begin{equation*}
	A_\mrm m \leq 1.
\end{equation*}


\subsection{Non-negativity,  lower and upper bounds of $ h_\mrm m $}\label{subsec:pointwise-bounds-h-app}

Let $ \underline h, \overline h \in [0,\infty) $ be the lower and upper bounds of $ h_\mrm{in} $, respectively, i.e., $ - 0 \leq \underline h \leq h_\mrm{in} \leq \overline h <\infty $ (see \eqref{cnts:initial-2}). In this section, we will show that
$$  \dfrac{1}{4} \underline h \leq h_\mrm{m} \leq 4 \overline h $$
locally in time. Again we assume that $ \vec u^o $ has the regularity \eqref{rg:24-feb-3}.

\vspace{5mm}

{\noindent\bf Non-negativity of $ h_\mrm m $: }

After applying the $ L^2 $-inner product of \eqref{eq:h-lin} with $ \bigl(- h_\mrm{m}^-\bigr) $ and applying integration by parts in the resultant, one has
\begin{equation}\label{eq:positivity-h-app}
	\dfrac{1}{2} \dfrac{d}{dt} \norm{h_\mrm{m}^-}{\Lnorm{2}}^2 = - \dfrac{1}{2} \int \abs{h_\mrm{m}^-}{2} \dv \vec u^o  \idx \lesssim \norm{\dv \vec u^o}{\Lnorm{\infty}} \norm{h_\mrm{m}^-}{\Lnorm{2}}^2,
\end{equation}
since the term $ \mathcal S_{h_\mrm m, \omega,\nu}(- h_\mrm{m}^- ) $ vanishes.
Therefore, applying Gr\"onwall's inequality to \eqref{eq:positivity-h-app}, as before in \eqref{eq:positivity-A-app}, eventually implies
\begin{equation*}
	h_\mrm{m} \geq 0. 
\end{equation*}

\vspace{5mm}

{\noindent\bf Lower and upper bounds of  $ h_\mrm{m} $: }

Since $ A_\mrm m  \in [0,1] $, one has $ \abs{\mathcal S_{h_\mrm m,\omega,\nu}}{} \leq 3 (\abs{\overline f}{} + \abs{\underline f}{})  $. Then following the characteristic method, 
since $ h_\mrm m \geq 0 $,
one has
\begin{gather*}
	\dt ( e^{-\int_0^t \norm{\dv \vec u^o}{\Lnorm{\infty}}(s) \,ds} h_\mrm m )+ \vec u^o \cdot \nabla ( e^{-\int_0^t \norm{\dv \vec u^o}{\Lnorm{\infty}}(s) \,ds} h_\mrm m)\\
	 \leq  3 (\abs{\overline f}{} + \abs{\underline f}{}) e^{-\int_0^t \norm{\dv \vec u^o}{\Lnorm{\infty}}(s) \,ds}\, .
\end{gather*}
Thus, integrating in the above inequation along the characteristic path given by $ \vec u^o $ yields
\begin{equation}\label{eq:upper-h-app}
h_\mrm m(\vec x, t) \leq \biggl( \overline h + 3 ( \abs{\overline f}{} + \abs{\underline f}{}) t \biggr) \times e^{\int_0^t \norm{\dv \vec u^o(s)}{\Lnorm{\infty}} \,ds}\, .
\end{equation}
Similarly, one can show that
\begin{equation}\label{eq:lower-h-app}
	h_\mrm m(\vec x, t) \geq \biggl( \underline h - 3 ( \abs{\overline f}{} + \abs{\underline f}{}) t \biggr) \times e^{-\int_0^t \norm{\dv \vec u^o(s)}{\Lnorm{\infty}} \,ds}\, .
\end{equation}
Then it immediately follows that $ \dfrac{1}{4} \underline h \leq h_\mrm m \leq 4 \overline h $ provided that the following conditions are satisfied:
\begin{equation}\label{cnd:small-time-001}
	\begin{gathered}
		0 < t \leq \begin{cases} 
		\dfrac{\underline h}{6(\abs{\overline f}{} + \abs{\underline f}{})} 
		& \text{if} ~ \underline h > 0, \\
		\\
		\dfrac{\overline h}{3(\abs{\overline f}{} + \abs{\underline f}{})} & \text{if} ~ \underline h = 0,
		 \end{cases}
		\\
		\text{and} ~~~~ ~~~~
		e^{\int_0^t \norm{\dv \vec u^o(s)}{\Lnorm{\infty}} \,ds} \leq  e^{t^{1/2} \bigl(\int_0^t \norm{\vec u^o(s)}{\Hnorm{3}}^2 \,ds\bigr)^{1/2}} \leq 2.
	\end{gathered}
\end{equation}

\subsection{Non-vanishing total ice mass}

Due to the fact that $ \abs{\mathcal S_{h_\mrm m,\omega,\nu}}{} \leq 3 (\abs{\overline f}{} + \abs{\underline f}{})  $, one can show immediately after integrating \eqref{eq:h-lin}, that
\begin{equation*}
	\dfrac{d}{dt} \int h_\mrm m \idx \leq 3 (\abs{\overline f}{} + \abs{\underline f}{}) \abs{\Omega}{}.
\end{equation*}
Therefore,
\begin{equation}\label{non-vanishi-mass}
	\dfrac{1}{2} \int h_\mrm{in} \idx \leq \int h_\mrm m \idx \leq 2 \int h_\mrm{in} \idx,
\end{equation}
provided
\begin{equation}\label{cnd:small-time-002}
	t \leq \dfrac{\int h_\mrm{in} \idx}{6(\abs{\overline f}{} + \abs{\underline f}{}) \abs{\Omega}{}}.
\end{equation}

\subsection{Well-posedness of \eqref{eq:u-lin} with strictly positive ice mass}\label{subsec:well-posedness_2.1c}

Consider $ h_\mrm{in} \geq \underline h > 0 $. Then we have shown  in section \ref{subsec:pointwise-bounds-h-app} that $ h_\mrm m \geq \underline h/4 > 0 $ locally in time. Then during this local time, \eqref{eq:u-lin} is a non-degenerate biharmonic evolutionary equation. Then following the standard Galerkin method, one can establish the well-posedness of strong solutions to \eqref{eq:u-lin}, provided that $ \vec u^o $ is sufficiently smooth. We omit the details here and refer interesting readers to \cite[section 7]{Evans}.

\section{Well-posedness of solutions to \eqref{eq:sea-ice-app} with $\underline{h} > 0$ and $ \iota > 0 $ fixed}\label{sec:first-layer-app}

In this section, we aim at showing that the map defined by
\begin{equation}\label{def:mapping}
	\mathfrak M: \vec u^o \mapsto \vec u_\mrm m,
\end{equation}
where $ \vec u_\mrm{m} $ is the unique solution to \eqref{eq:u-lin} with $ h_\mrm m $ and $ A_\mrm m $ being solutions to \eqref{eq:h-lin} and \eqref{eq:A-lin}, respectively, is bounded in $ \mathfrak X_{T^*} $ and contracting with contraction constant $ 1/2 $ in $ L^\infty(0,T^*;L^{2}(\Omega))\cap L^2(0,T^*;H^{2}(\Omega)) $, where

\begin{equation}\label{def:compact-space}
	\begin{gathered}
	\mathfrak X_{T^*} : = \bigl\lbrace \vec u| \vec u \in L^\infty(0,T^*;\Hnorm{2}) \cap L^2(0,T^*;\Hnorm{3}), \\
	 \dt \vec u, \nabla^4 \vec u \in L^2(0,T^*;\Lnorm{2}) \bigr\rbrace,
	\end{gathered}
\end{equation}
for some $ T^* \in (0,\infty) $ to be determined. Throughout this section, unless stated otherwise, the initial data for $ \vec u_\mrm m, \vec u^o, h $, and $ A $ are assumed to be $ \vec u_\mrm{in}, \vec u_\mrm{in},h_\mrm{in} $, and $ A_\mrm{in} $, given in Theorem \ref{thm:well-posedness-app-2}, respectively.

Consequently, one can apply the Banach fixed-point theorem, i.e., the contraction mapping theorem, to show the existence of solutions to system \eqref{eq:sea-ice-app}.


Let $ \mathfrak c_{\mrm{in}} \in (0,\infty) $ be the bound of the initial data defined by
\begin{equation}\label{def:initial-bound}
	\norm{\nabla h_\mrm{in}, \nabla A_\mrm{in}}{\Lnorm{4}} + \norm{\vec u_{\mrm {in}}}{\Hnorm{2}} \leq \mathfrak c_{\mrm{in}}
\end{equation}

\subsection{Uniform bounds}\label{subsec:uniform-bound}

Let $ \vec u^o \in \mathfrak X_{T^*} $ satisfy
\begin{equation}\label{def:solution-bound}
	\sup_{0\leq s \leq t}\norm{\vec u^o(s)}{\Hnorm{2}}^2 + \int_0^t \bigl( \norm{\dt \vec u^o(s)}{\Lnorm{2}}^2 + \norm{\vec u^o(s)}{\Hnorm{3}}^2 \bigr) \,ds \leq \mathfrak c_o,
\end{equation}
with $ t \in [0,T^*] $,
for some $ \mathfrak c_o \in (0,\infty) $ to be determined later.

\vspace{5mm}

{\noindent\bf Estimates for $ h_\mrm m $ and $ A_\mrm m $}
Aside from the point-wise estimates deduced in Sections \ref{subsec:pointwise-A-app} and \ref{subsec:pointwise-bounds-h-app}, we shall need a uniform $ H^1 $-estimate for $ A_\mrm m $ and $ h_\mrm m $.

We record the equation after applying 
$ \partial \in \lbrace \partial_x, \partial_y \rbrace $ to \eqref{eq:h-lin}, as follows:
\begin{equation}\label{eq:d-h-lin}
	\dt \partial h_\mrm m + \vec u^o \cdot  \nabla \partial h_\mrm m + \partial \vec u^o \cdot \nabla h_\mrm m + \partial h_\mrm m \dv \vec u^o + h_\mrm m \dv \partial \vec u^o = \partial \mathcal S_{h_\mrm m, \omega,\nu}.
\end{equation}
Then taking the $ L^2 $-inner product of \eqref{eq:d-h-lin} with $ 4 \abs{\partial h_\mrm m}{2} \partial h_\mrm m $ leads to, after applying integration by parts,
\begin{equation}\label{est:h-lin-001}
\begin{aligned}
	& \dfrac{d}{dt}\norm{\partial h_\mrm m}{\Lnorm{4}}^4 = - 3 \int \dv \vec u^o \abs{\partial h_\mrm m}{4} \idx \\
	& ~~ - 4 \int \bigl( \partial \vec u^o \cdot \nabla h_\mrm m + h_\mrm m \dv \partial \vec u^o \bigr) \abs{\partial h_\mrm m}{2} \partial h_\mrm m \idx + 4 \int \partial \mathcal S_{h_\mrm m, \omega,\nu} \abs{\partial h_\mrm m}{2} \partial h_\mrm m \idx \\
	& \lesssim  \norm{\nabla \vec u^o}{\Lnorm{\infty}} \norm{\nabla h_\mrm m}{\Lnorm{4}}^4 + \norm{h_\mrm m}{\Lnorm{\infty}} \norm{\nabla^2 \vec u^o}{\Lnorm{4}} \norm{\nabla h_\mrm m}{\Lnorm{4}}^3\\
	& ~~~~ + \int \abs{\partial \mathcal S_{h_\mrm m, \omega,\nu}}{} \abs{\partial h_\mrm m}{2} \partial h_\mrm m \idx  .
\end{aligned}
\end{equation}
Meanwhile, simple calculation shows that
\begin{equation*}
	\abs{\partial \mathcal S_{h_\mrm m, \omega,\nu}}{} \lesssim \bigl( \dfrac{1}{\omega} + \dfrac{1}{\nu^{1/2}} \bigr) \abs{\partial h_\mrm m}{} + \bigl( 1 + \dfrac{\abs{h_\mrm m }{}}{\omega^2} \bigr) \abs{\partial A_\mrm m}{},
\end{equation*}
where we have used \eqref{cnts:f}. Consequently, one concludes from \eqref{est:h-lin-001} that
\begin{equation}\label{est:h-lin-002}
	\begin{aligned}
	& \dfrac{d}{dt} \norm{\nabla h_\mrm m}{\Lnorm{4}}^4 \lesssim \bigl( \norm{\nabla \vec u^o}{\Lnorm{\infty}} + \dfrac{1}{\omega} + \dfrac{1}{\nu} \bigr) \norm{\nabla h_\mrm m}{\Lnorm{4}}^4 \\
	& ~~~~ + \bigl(1 + \dfrac{\norm{h_\mrm m}{\Lnorm{\infty}}}{\omega^2} \bigr) \norm{\nabla A_\mrm m}{\Lnorm{4}}\norm{\nabla h_\mrm m}{\Lnorm{4}}^3 \\
	& ~~~~+ \norm{h_\mrm m}{\Lnorm{\infty}} \norm{\nabla^2 \vec u^o}{\Lnorm{4}} \norm{\nabla h_\mrm m}{\Lnorm{4}}^3.
	\end{aligned}
\end{equation}


The estimate for $ \nabla A_\mrm m $ is obtained from \eqref{eq:A-lin} in a similar fashion, we record it here:
\begin{equation}\label{est:A-lin-001}
	\begin{aligned}
		& \dfrac{d}{dt} \norm{\nabla A_\mrm m}{\Lnorm{4}}^4 \lesssim \biggl( \norm{\nabla \vec u^o}{\Lnorm{\infty}}+ \dfrac{\norm{\nabla \vec u^o}{\Lnorm{\infty}}}{\omega} \\
		& ~~~~ ~~~~ + 1 + \dfrac{1}{\nu} + \dfrac{\norm{h_\mrm m}{\Lnorm{\infty}}}{\omega^2}  \biggr)\norm{\nabla A_\mrm m}{\Lnorm{4}}^4 \\
		& ~~~~ +  \bigl( \dfrac{1}{\nu^2} + \dfrac{1}{\omega^2} \bigr) \norm{\nabla h_\mrm m}{\Lnorm{4}} \norm{\nabla A_\mrm m}{\Lnorm{4}}^3
		\\ & ~~~~
		 + \norm{\nabla^2 \vec u^o}{\Lnorm{4}} \norm{\nabla A_\mrm m}{\Lnorm{4}}^3,
	\end{aligned}
\end{equation}
where we have used the fact that $ A_\mrm m \in [0,1] $.

After combining \eqref{est:h-lin-002} and \eqref{est:A-lin-001}
and
applying Gr\"onwall's inequality,  one can derive that
\begin{equation}\label{est:h-A-lin-h1-total}
	\sup_{0\leq s \leq t}\norm{\nabla h_\mrm m(s), \nabla A_\mrm m(s)}{\Lnorm{4}}^4 \leq e^{H_{h,A,1}(t)} \bigl( \norm{\nabla h_\mrm{in}, \nabla A_\mrm{in}}{\Lnorm{4}}^4 + G_{h,A,1}(t)  \bigr),
\end{equation}
where
\begin{align}
	& \begin{aligned}\label{est:h-A-lin-h1-1}
		& H_{h,A,1}(t) : = C_{\omega, \nu}\int_0^t \bigl( 1 + \norm{\nabla \vec u^o(s)}{\Lnorm{\infty}} + \norm{h_\mrm m(s)}{\Lnorm{\infty}} \\
		& ~~~~ + \norm{\nabla^2 \vec u^o(s)}{\Lnorm{4}} + \norm{h_\mrm m(s)}{\Lnorm{\infty}} \norm{\nabla^2 \vec u^o(s)}{\Lnorm{4}} \bigr) \,ds,
	\end{aligned}\\
	&
	\label{est:h-A-lin-h1-2}
		G_{h,A,1}(t) : = \int_0^t \bigl( 1 +  \norm{h_\mrm m(s)}{\Lnorm{\infty}} + \norm{\nabla^2 \vec u^o(s)}{\Lnorm{4}}  \bigr) \,ds.
\end{align}

On the other hand, in direct consequence of equations \eqref{eq:h-lin} and \eqref{eq:A-lin}, one has
\begin{equation}\label{est:h-A-lin-dt}
	\begin{gathered}
		\norm{\dt h_\mrm m, \dt A_\mrm m}{\Lnorm{4}} \leq C \bigl( 1 + 1/\nu + \norm{\nabla \vec u^o}{\Lnorm{4}} \\
		~~~~ + \norm{h_\mrm m}{\Lnorm{\infty}} \norm{\nabla \vec u^o}{\Lnorm{4}}
		 + \norm{\vec u^o}{\Lnorm{\infty}}\norm{\nabla h_\mrm m, \nabla A_\mrm m}{\Lnorm{4}} \bigr),
	\end{gathered}	
\end{equation}
where we have used the fact that $ 0\leq A_\mrm m \leq 1 $ and \eqref{cnts:f}.

\vspace{5mm}
{\noindent\bf Estimates for $ \vec u_\mrm m $}

Taking the $ L^2 $-inner product of \eqref{eq:u-lin} with $ 2 \vec u_\mrm m + 2 \dt \vec u_\mrm m  - 2 \Delta \vec u_{\mrm m} $ leads to, after applying integration by parts,
\begin{equation}\label{est:u-lin-001}
	\begin{aligned}
		& \dfrac{d}{dt} \norm{\rho_\mrm{ice}^{1/2}h_\mrm m ^{1/2} \vec u_\mrm m, \iota^{1/2} \nabla^2 \vec u_\mrm m, \rho_\mrm{ice}^{1/2}h^{1/2}_\mrm m \nabla \vec u_\mrm m }{\Lnorm{2}}^2 \\
		& ~~~~ ~~~~ + 2 \norm{\rho_\mrm{ice}^{1/2}h_\mrm m ^{1/2} \dt \vec u_\mrm m, \iota^{1/2} \nabla^2 \vec u_\mrm m, \iota^{1/2} \nabla^3 \vec u_\mrm m }{\Lnorm{2}}^2 \\
		&  = \underbrace{\int \rho_\mrm{ice} \dt h_\mrm m \abs{\vec u_\mrm m}{2} \idx}_{\mathcal R_{1}} + \underbrace{2 \int \rho_\mrm{ice} (\nabla h_\mrm m \cdot \nabla)\vec  u_\mrm m \cdot \dt \vec u_\mrm m\idx}_{\mathcal R_{2}} \\
		& ~~ \underbrace{- \int \rho_\mrm{ice} \dt h_\mrm m \abs{\nabla \vec u_\mrm m}{2}\idx}_{\mathcal R_{3}}  \underbrace{- 2 \int \rho_\mrm{ice} h_\mrm m (\vec u^o \cdot \nabla) \vec u^o \cdot (\vec u_\mrm m + \dt \vec u_\mrm m - \Delta \vec u_\mrm m) \idx}_{\mathcal R_{4}}  \\
		& ~~ \underbrace{- 2 \int \nabla p_\mrm m \cdot (\vec u_\mrm m + \dt \vec u_\mrm m- \Delta \vec u_\mrm m) \idx}_{\mathcal R_{5}} + \underbrace{2 \int \mathcal F_\mrm m \cdot(\vec u_\mrm m + \dt \vec u_\mrm m- \Delta \vec u_\mrm m) \idx}_{\mathcal R_{6}} \\
		& ~~ + \underbrace{2 \int \dv \mathbb S_{\varepsilon, \mu,\lambda,\mrm m} \cdot (\vec u_\mrm m + \dt \vec u_\mrm m- \Delta \vec u_\mrm m) \idx}_{\mathcal R_{7}}. 
	\end{aligned}
\end{equation}
We obtain the following estimates for the $ \mathcal R_j $ terms by
applying H\"older's inequality and the Sobolev embedding inequality:
\begin{align*}
	& \mathcal R_1 \lesssim \norm{\dt h_\mrm m}{\Lnorm 2} \norm{\vec u_{\mrm m}}{\Lnorm{4}}^2, \\
	& \mathcal R_2 \lesssim \norm{\dt \vec u_\mrm m}{\Lnorm{2}} \norm{\nabla \vec u_\mrm m}{\Lnorm{4}} \norm{\nabla h_\mrm m}{\Lnorm{4}} , \\
	& \mathcal R_3 \lesssim \norm{\dt h_\mrm m}{\Lnorm{2}} \norm{\nabla \vec u_\mrm m}{\Lnorm{4}}^2, \\
	& \mathcal R_4 \lesssim \norm{h_\mrm m}{\Lnorm{\infty}} \norm{\vec u^o}{\Lnorm{4}} \norm{\nabla \vec u^o}{\Lnorm{4}} \norm{\dt \vec u_\mrm m, \vec u_\mrm m, \nabla^2 \vec u_\mrm m}{\Lnorm{2}}, \\
	& \mathcal R_5 \lesssim  \norm{\nabla p_\mrm m}{\Lnorm{2}}\norm{\dt \vec u_\mrm m, \vec u_\mrm m, \nabla^2 \vec u_\mrm m}{\Lnorm{2}}, \\
	& \mathcal R_6 \lesssim \bigl(1 + \norm{\vec u^o}{\Lnorm{2}} + \norm{h_\mrm m}{\Lnorm{\infty}} \norm{\vec u^o}{\Lnorm{2}}\bigr) \norm{\dt \vec u_\mrm m, \vec u_\mrm m, \nabla^2 \vec u_\mrm m}{\Lnorm{2}}, \\
	& \mathcal R_7 \lesssim \bigl( \dfrac{1}{\varepsilon} \norm{p_\mrm m}{\Lnorm{\infty}} \norm{\nabla^2 \vec u^o}{\Lnorm{2}} + (\mu +\lambda) \norm{\nabla^2 \vec u^o}{\Lnorm{2}} + \norm{\nabla p_\mrm m}{\Lnorm{2}}
	 \bigr) \\
	 & \qquad\qquad\qquad \times \norm{\dt \vec u_\mrm m, \vec u_\mrm m, \nabla^2 \vec u_\mrm m}{\Lnorm{2}}.
\end{align*}

To deduce the above estimates, consider $ \underline h > 0 $ and let $ t $ satisfy \eqref{cnd:small-time-001} and \eqref{cnd:small-time-002}. Therefore, the estimates in Section \ref{subsec:pointwise-bounds-h-app} guarantee that $ 0 < 1/4 \underline h \leq h_\mrm m \leq 4 \overline h < \infty $.
Consequently, \eqref{est:u-lin-001} yields, after applying the Sobolve embedding inequality and H\"older's inequality,
\begin{equation}\label{est:u-lin-002}
	\begin{aligned}
		& \dfrac{d}{dt} \norm{\rho_\mrm{ice}^{1/2}h_\mrm m ^{1/2} \vec u_\mrm m, \rho_\mrm{ice}^{1/2}h^{1/2}_\mrm m \nabla \vec u_\mrm m , \iota^{1/2} \nabla^2 \vec u_\mrm m}{\Lnorm{2}}^2 \\
		& ~~~~ ~~~~ + 2 \norm{\rho_\mrm{ice}^{1/2}h_\mrm m ^{1/2} \dt \vec u_\mrm m, \iota^{1/2} \nabla^2 \vec u_\mrm m, \iota^{1/2} \nabla^3 \vec u_\mrm m }{\Lnorm{2}}^2\\
		& \leq C_{\varepsilon,\mu,\lambda,\underline h,\overline h} \bigl( \norm{\dt h_\mrm m}{\Lnorm{4}} + \norm{\nabla h_\mrm m,\nabla A_\mrm m}{\Lnorm{4}}^2  \bigr)\\
		& ~~~~ \times \bigl( \norm{\rho_\mrm{ice}^{1/2}h_\mrm m ^{1/2} \vec u_\mrm m, \rho_\mrm{ice}^{1/2}h^{1/2}_\mrm m \nabla \vec u_\mrm m, \iota^{1/2} \nabla^2 \vec u_\mrm m}{\Lnorm{2}}^4 + 1 \bigr).
	\end{aligned}
\end{equation}
Furthermore, consider $ t $ small enough such that 
\begin{equation}\label{cnd:small-time-003}
	\begin{aligned}
		&
		H_{h,A,1}(t) + G_{h,A,1}(t)
		\leq C_{\omega, \nu, \overline h} t^{1/2} \bigl( t^{1/2} + \bigl( \int_0^t \norm{\nabla \vec u^o(s)}{\Hnorm{2}}^2 \,ds\bigr)^{1/2}\bigr) \\
		& ~~~~ ~~~~ ~~~~ \leq C_{\omega, \nu, \overline h} t^{1/2} \bigl( t^{1/2} + \mathfrak c_o^{1/2} \bigr)   \leq 1,
	\end{aligned}
\end{equation}
where we have applied H\"older's inequality.
Then \eqref{est:h-A-lin-h1-total} and \eqref{est:h-A-lin-dt} imply that, after applying the Sobolev embedding inequality,
\begin{equation}\label{est:h-A-lin-total}
	\norm{\nabla h_\mrm m, \nabla A_\mrm m,\dt h_\mrm m, \dt A_\mrm m}{\Lnorm{4}} \leq C_{\omega,\nu, \overline h,\mathfrak c_\mrm{in}}\bigl( 1 + \mathfrak c_o^{1/2} \bigr).
\end{equation}
Consequently, 
\eqref{est:u-lin-002} yields the following estimate:
\begin{equation}\label{est:u-lin-total}
	\begin{aligned}
		& \sup_{0\leq s\leq t}\norm{\vec u_\mrm m(s)}{\Hnorm{2}}^2 + \int_0^t \bigl( \norm{\dt \vec u_\mrm m(s)}{\Lnorm{2}}^2 + \norm{\vec u_\mrm m(s)}{\Hnorm{3}}^2 \bigr) \,ds \\
		&  \leq C_{\varepsilon, \iota, \mu, \lambda,\omega, \nu, \underline h ,\overline h, \mathfrak c_\mrm{in}} \biggl\lbrack \bigl( \dfrac{C_{\varepsilon, \mu, \lambda,\omega, \nu, \underline h ,\overline h, \mathfrak c_\mrm{in},1}}{C_{\varepsilon, \mu, \lambda,\omega, \nu, \underline h ,\overline h, \mathfrak c_\mrm{in},2} - (1+\mathfrak c_o)t} - 1 \bigr)^2 \bigl(1 + (1+\mathfrak c_o) t \bigr) +1 \biggr\rbrack\\
		& \leq C_{\varepsilon, \iota, \mu, \lambda,\omega, \nu, \underline h ,\overline h, \mathfrak c_\mrm{in}} \biggl\lbrack 2 \bigl( 2 \dfrac{C_{\varepsilon, \mu, \lambda,\omega, \nu, \underline h ,\overline h, \mathfrak c_\mrm{in},1}}{C_{\varepsilon, \mu, \lambda,\omega, \nu, \underline h ,\overline h, \mathfrak c_\mrm{in},2}} - 1\bigr)^2+ 1 \biggr\rbrack,
	\end{aligned}
\end{equation}
provided that $ t $ is small enough and where we have made the choice
\begin{equation}\label{def:uniform-u-lin}
	 \mathfrak c_o : = C_{\varepsilon, \iota, \mu, \lambda,\omega, \nu, \underline h ,\overline h, \mathfrak c_\mrm{in}} \biggl\lbrack 2 \bigl( 2 \dfrac{C_{\varepsilon, \mu, \lambda,\omega, \nu, \underline h ,\overline h, \mathfrak c_\mrm{in},1}}{C_{\varepsilon, \mu, \lambda,\omega, \nu, \underline h ,\overline h, \mathfrak c_\mrm{in},2}} - 1\bigr)^2+ 1 \biggr\rbrack,
\end{equation}
where the right-hand side is as in \eqref{est:u-lin-total}. Then \eqref{cnd:small-time-001}, \eqref{cnd:small-time-002}, \eqref{cnd:small-time-003}, and \eqref{est:u-lin-total} imply that,
there exists $ T^* \in (0,\infty) $ such that
\begin{subequations}
\begin{equation}\label{est:u-lin-bound}
	 \sup_{0\leq s\leq t}\norm{\vec u_\mrm m(s)}{\Hnorm{2}}^2 + \int_0^t \bigl( \norm{\dt \vec u_\mrm m(s)}{\Lnorm{2}}^2 + \norm{\vec u_\mrm m(s)}{\Hnorm{3}}^2 \bigr) \,ds \leq \mathfrak c_o,
\end{equation}
and
\begin{equation}
	\dfrac{1}{4} \underline h \leq h_\mrm m \leq 4 \overline h, ~~ \dfrac{1}{2} \int h_\mrm{in} \idx \leq \int h_\mrm m \idx \leq 2 \int h_\mrm{in} \idx,
\end{equation}
\end{subequations}
for $ t \in [0,T^*] $. In addition, using equation \eqref{eq:sea-ice-app}, it is easy to obtain
\begin{equation}\label{est:u-lin-bounds-2}
	\int_0^t \norm{\Delta^2 \vec u_\mrm{m}(s)}{\Lnorm{2}}^2 \,ds \leq C_{\overline h, \mu, \lambda, \varepsilon} \mathfrak c_o.
\end{equation}
Therefore, $ \mathfrak M $, defined in \eqref{def:mapping}, maps $ \mathfrak X_{T^*} $ into itself for such choices of $ T^* $ and $ \mathfrak c_o $.

We remark here that, $ \mathfrak c_0 \rightarrow \infty $ as $ \iota \rightarrow 0^+ $, i.e., the estimates we obtain here depend on $ \iota >0 $. We will remove the dependency of $ \iota $ in Section \ref{sec:well-posedness}.

\subsection{Contraction mapping and well-posedness}\label{subsec:contracting-mapping}

For $ j = 1 ,2 $,
consider $ \vec u^o_j  \in \mathfrak X_{T^*} $ satisfying \eqref{def:solution-bound}, and let $ h_{\mrm m,j}, A_{\mrm m, j} $, and $ \vec u_{\mrm m, j} = \mathfrak M (\vec u^o_j) $, be the solutions to \eqref{eq:h-lin}, \eqref{eq:A-lin}, and \eqref{eq:u-lin}, respectively, with $ \vec u^o $ replaced by $ \vec u^o_j $ and with the same initial data.
Then we have the estimates of $ h_{\mrm{m}, j}, A_{\mrm m, j} $, and $ \vec u_{\mrm m, j} $ as in Sections \ref{subsec:pointwise-A-app} and \ref{subsec:pointwise-bounds-h-app}, as well as \eqref{est:h-A-lin-total} and \eqref{est:u-lin-bound}.

In the following, let $ \sigma \in (0,1) $ be a constant to be determined later.
Denote by
\begin{equation}\label{def:delta-comtracting}
	\begin{gathered}
	\delta h_\mrm m := h_{\mrm m,1} - h_{\mrm m,2}, \quad \delta A_\mrm m := A_{\mrm m,1} - A_{\mrm m,2},\\
	\delta \vec u_\mrm m := \vec u_{\mrm m,1} - \vec u_{\mrm m,2},\quad \delta \vec u^o := \vec u^o_{1} - \vec u^o_{2}.
	\end{gathered}
\end{equation}
The notations
\begin{gather*}
\delta p_\mrm m,
\delta \mathbb S_{\varepsilon,\mu, \lambda, \mrm m}, \delta\mathcal F_\mrm m, \delta \mathcal S_{h_\mrm m, \omega, \nu}, \delta \mathcal S_{A_\mrm m, \omega, \nu}, \delta \chi_{A_\mrm m}^\omega,
\end{gather*}
have similar meanings.
Then $ \delta h_\mrm m, \delta A_\mrm m, \delta \vec u_\mrm m $ satisfy
\begin{subequations}
\begin{gather}
	\dt \delta h_\mrm m + \dv  ( \delta h_\mrm m \vec u^o_1)  + \dv  ( h_{\mrm{m},2} \delta \vec u^o ) = \delta \mathcal S_{h_\mrm m, \omega,\nu},
	\label{eq:h-diff-lin}
	\\
	\begin{gathered}
	\dt \delta A_\mrm m + \dv  ( \delta A_\mrm m \vec u^o_1)  + \dv  ( A_{\mrm{m},2} \delta \vec u^o ) = \delta \mathcal S_{A_\mrm m, \omega,\nu} \\
	 + \delta A_{\mrm m} \dv \vec u^o_1\cdot \chi_{A_{\mrm m, 1}}^\omega
	 +  A_{\mrm m,2} \dv \delta  \vec u^o\cdot \chi_{A_{\mrm m, 1}}^\omega + A_{\mrm m,2} \dv   \vec u^o_2 \cdot \delta \chi_{A_{\mrm m}}^\omega,
	 \end{gathered}
	 \label{eq:A-diff-lin}
	 \\
	 \begin{gathered}
	 \rho_{\mrm{ice}}h_{\mrm m, 1} \dt  \delta \vec u_\mrm m  + \rho_\mrm{ice} \delta h_\mrm m \dt \vec u_{\mrm m,2} + \iota \Delta^2 \delta \vec u_\mrm m = - \rho_\mrm{ice} h_{\mrm m,1} \vec u^o_1 \cdot \nabla \delta \vec u^o \\
	  - \rho_\mrm{ice} h_{\mrm m,1} \delta \vec u^o \cdot \nabla  \vec u^o_2 - \rho_\mrm{ice}\delta h_{\mrm m} \vec u^o_2 \cdot \nabla  \vec u^o_2 - \nabla \delta p_\mrm m + \dv \delta \mathbb S_{\varepsilon,\mu,\lambda,\mrm m} + \delta \mathcal F_\mrm m.
	  \end{gathered}
	  \label{eq:u-diff-lin}
\end{gather}
\end{subequations}

After taking the $ L^2 $-inner product of \eqref{eq:h-diff-lin} and \eqref{eq:A-diff-lin} with $ 4 \abs{\delta h_\mrm m}{2}  \delta h_{\mrm m} $ and $ 4 \abs{\delta A_\mrm m}{2}\delta A_\mrm m $, respectively, and applying integration by parts in the resultant, one has
\begin{equation}\label{est:24-feb-1}
\begin{aligned}
	& \dfrac{d}{dt} \norm{\delta h_\mrm m,\delta A_\mrm m}{\Lnorm{4}}^4 = \underbrace{- 3 \int \dv \vec u_1^o ( \abs{\delta h_\mrm m}{4} +  \abs{\delta A_\mrm m}{4}) \idx}_{\mathcal R_{8}} \\
	& ~~~~ \underbrace{ - 4 \int ( \delta \vec u^o \cdot \nabla h_{\mrm m, 2} \abs{\delta h_\mrm m}{2}  \delta h_{\mrm m} + \delta \vec u^o \cdot \nabla A_{\mrm m, 2} \abs{\delta A_\mrm m}{2}  \delta A_{\mrm m})\idx}_{\mathcal R_{9}} \\
	& ~~~~ \underbrace{- 4 \int ( h_{\mrm m,2} \dv \delta \vec u^o \abs{\delta h_\mrm m}{2}  \delta h_{\mrm m} + A_{\mrm m,2} \dv \delta \vec u^o \abs{\delta A_\mrm m}{2}  \delta A_{\mrm m})\idx}_{\mathcal R_{10}} \\
	& ~~~~ + \underbrace{4 \int \dv \vec u_1^o \abs{\delta A_\mrm m}{4} \chi_{A_{\mrm m,1}}^\omega \idx}_{\mathcal R_{11}} + \underbrace{4 \int A_{\mrm m,2} \dv \delta \vec u^o \abs{\delta A_\mrm m}{2} \delta A_\mrm m \chi_{A_{\mrm m,1}}^\omega \idx}_{\mathcal R_{12}} \\
	& ~~~~ +\underbrace{4 \int A_{\mrm m,2} \dv\vec u_2^o \abs{\delta A_\mrm m}{2} \delta A_\mrm m \delta \chi_{A_{\mrm m}}^\omega \idx}_{\mathcal R_{13}} + \underbrace{4 \int \delta \mathcal S_{h_\mrm m,\omega,\nu}\abs{\delta h_\mrm m}{2}  \delta h_{\mrm m}\idx}_{\mathcal R_{14}}  \\
	& ~~~~ + \underbrace{4 \int  \delta \mathcal S_{A_\mrm m,\omega,\nu}\abs{\delta A_\mrm m}{2}  \delta A_{\mrm m}\idx}_{\mathcal R_{15}}.
\end{aligned}
\end{equation}
In the following, we sketch the estimates of the $ \mathcal R_j $ terms by applying H\"older's inequality and the Sobolev embedding inequality:
\begin{equation}\label{est:24-feb-2}
\begin{aligned}
	& \mathcal R_8 + \mathcal R_{11} + \mathcal R_{13} \lesssim \biggl(\norm{\dv \vec u^o_1}{\Lnorm{\infty}} + (\dfrac{1}{\omega}+ \dfrac{1}{\omega^{2}})\norm{\dv \vec u^o_2}{\Lnorm{\infty}}\biggr)\\
	& \qquad\qquad \times \norm{\delta h_\mrm m, \delta A_\mrm m}{\Lnorm{4}}^4,\\
	& \mathcal R_9 \lesssim \norm{\delta \vec u^o}{\Lnorm{\infty}} \norm{\nabla h_{\mrm m,2},\nabla A_{\mrm m,2}}{\Lnorm{4}} \norm{\delta h_\mrm m, \delta A_\mrm m}{\Lnorm{4}}^3,\\
	& \mathcal R_{10} + \mathcal R_{12} \lesssim (\overline h + 1) \norm{\dv \delta \vec u^o}{\Lnorm{4}} \norm{\delta h_\mrm m, \delta A_\mrm m}{\Lnorm{4}}^3,\\
	& \mathcal R_{14} + \mathcal R_{15} \lesssim C_{\overline h,\omega,\nu} \norm{\delta h_\mrm m, \delta A_\mrm m}{\Lnorm{4}}^4,
\end{aligned}
\end{equation}
where we have used the identity
\begin{equation*}
	\delta\biggl( \dfrac{g}{g+ \varepsilon} \biggr) = \dfrac{\delta g}{g_1+ \varepsilon}-\dfrac{g_2 \delta g}{(g_1+\varepsilon)(g_2+\varepsilon)}
\end{equation*}
for $ g = (1-A_\mrm{m})^+ = 1 - A_\mrm{m} $ in the estimate of $ \delta \chi_{A_\mrm m}^\omega $ in $\mathcal R_{13}$.
In view of \eqref{est:24-feb-1} and \eqref{est:24-feb-2},
one has
\begin{equation}\label{est:h-A-diff-lin-001}
	\begin{gathered}
		\dfrac{d}{dt}\norm{\delta h_\mrm m,\delta A_\mrm m}{\Lnorm{4}}^2 \leq C_{\sigma,\varepsilon, \omega, \nu, \overline h, \mathfrak c_o,\mathfrak c_\mrm{in}}
		\norm{\delta h_\mrm m,\delta A_\mrm m}{\Lnorm{4}}^2
		+ \sigma \norm{\delta \vec u^o}{\Hnorm{2}}^2,
	\end{gathered}
\end{equation}
where we have used \eqref{def:solution-bound} and \eqref{est:h-A-lin-h1-total}.
Consequently, applying Gr\"onwall's inequality to \eqref{est:h-A-diff-lin-001} yields
\begin{equation}\label{est:h-A-diff-lin-002}
	\begin{gathered}
	\sup_{0\leq s\leq t} \norm{\delta h_\mrm m(s),\delta A_\mrm m(s)}{\Lnorm{4}}^2 \\  \leq \sigma \biggl( \int_0^t \norm{\delta \vec{u}^o(s)}{\Hnorm{2}}^2 \,ds \biggr)
	  e^{C_{\sigma,\varepsilon, \omega, \nu, \overline h, \mathfrak c_o,\mathfrak c_\mrm{in}} (t + t^{1/2})},
	 \end{gathered}
\end{equation}
where we have also employed Young's inequality.

Taking the $ L^2 $-inner product of \eqref{eq:u-diff-lin} with $ 2 \delta \vec u_\mrm m $ and applying integration by parts in the resultant yields
\begin{equation}\label{est:u-diff-lin-001}
	\begin{aligned}
		& \rho_\mrm{ice} \dfrac{d}{dt} \norm{h_{\mrm m,1}^{1/2} \delta \vec u_\mrm m}{\Lnorm{2}}^2 + 2 \iota \norm{ \nabla^2 \delta \vec u_\mrm m}{\Lnorm{2}}^2 = \underbrace{\rho_\mrm{ice} \int \dt h_{\mrm m,1} \abs{\delta \vec u_\mrm m }{2} \idx}_{\mathcal R_{16}}  \\
		& ~~~~ \underbrace{- 2 \int \rho_{\mrm{ice}} \delta h_\mrm m \dt \vec u_{\mrm{m},2} \cdot \delta \vec u_\mrm m \idx}_{\mathcal R_{17}} \underbrace{
		  - 2 \int \rho_\mrm{ice} h_{\mrm m,1} (\vec u_1^o \cdot \nabla ) \delta \vec u^o \cdot \delta \vec u_\mrm m \idx}_{\mathcal R_{18}}
		  \\ & ~~~~
		  \underbrace{- 2 \int \rho_\mrm{ice} h_{\mrm m, 1} (\delta \vec u^o \cdot \nabla) \vec u_2^o \cdot \delta \vec u_\mrm m \idx }_{\mathcal R_{19}}
		 \underbrace{ - 2 \int \rho_\mrm{ice} \delta h_\mrm m (\vec u_2^o \cdot \nabla) \vec u_2^o \cdot \delta \vec u_\mrm m \idx }_{\mathcal R_{20}}
		  \\ & ~~~~
		 + \underbrace{2 \int \delta p_\mrm m \dv \delta \vec u_\mrm m \idx }_{\mathcal R_{21}}
		 + \underbrace{2 \int \dv \delta \mathbb S_{\varepsilon,\mu,\lambda,\mrm m} \cdot \delta \vec u_\mrm m \idx }_{\mathcal R_{22}}
		 \\ & ~~~~
		 + \underbrace{2 \int \delta \mathcal F_\mrm m\cdot\delta \vec u_\mrm m \idx}_{\mathcal R_{23}}. 
	\end{aligned}
\end{equation}
In the following, again, we sketch the estimates for the terms $ \mathcal R_j $ by applying H\"older's inequality, the Gagliardo-Nirenberg inequality, and the Sobolev embedding inequality:
\begin{align*}
	& \mathcal R_{16} \lesssim \norm{\dt h_{\mrm m,1}}{\Lnorm{4}}\norm{\delta \vec u_\mrm m}{\Lnorm{2}}^{3/2}\norm{\delta \vec u_\mrm m}{\Hnorm{1}}^{1/2}
	,
	\\ 
	& \mathcal R_{17} \lesssim \norm{\dt \vec u_{\mrm m,2}}{\Lnorm{2}} \norm{\delta h_\mrm m}{\Lnorm{4}} \norm{\delta \vec u_\mrm m}{\Lnorm{2}}^{1/2} \norm{\delta \vec u_\mrm m}{\Hnorm{1}}^{1/2},
	\\ 
	& \mathcal R_{18} \lesssim \overline h \norm{\vec u_1^o}{\Hnorm{2}} \norm{\nabla \delta \vec u^o}{\Lnorm{2}} \norm{\delta \vec u_\mrm m}{\Lnorm{2}}  ,
	\\ 
	& \mathcal R_{19} \lesssim \overline h \norm{\delta \vec u^o}{\Lnorm{2}}^{1/2}\norm{\delta \vec u^o}{\Hnorm{1}}^{1/2} \norm{\nabla \vec u_2^o}{\Lnorm{4}}\norm{\delta \vec u_\mrm m}{\Lnorm{2}} ,
	\\ 
	& \mathcal R_{20} \lesssim \norm{\delta h_\mrm m}{\Lnorm{4}} \norm{\vec u_2^o}{\Hnorm{2}}\norm{\nabla \vec u_2^o}{\Lnorm{4}} \norm{\delta\vec u_\mrm m}{\Lnorm{2}} ,
	\\ 
	& \mathcal R_{21} \lesssim (1 + \overline h) \norm{\delta h_\mrm m, \delta A_\mrm m}{\Lnorm{2}} \norm{\nabla \delta \vec u_\mrm m}{\Lnorm{2}} ,
	\\
	& \mathcal R_{23} \lesssim  ( 1 + \overline h + \sum_{j=1}^2\norm{\vec u^o_j}{\Hnorm{2}} )(\norm{\delta \vec u^o}{\Lnorm{2}} + \norm{\delta h_\mrm m}{\Lnorm{2}}) \norm{\delta \vec u_\mrm m}{\Lnorm{2}}.
\end{align*}
To estimate $ \mathcal R_{22} $, we rewrite it as
\begin{align*}
	& \mathcal R_{22} = 2 \int \dv \bigl\lbrack \mu ( \nabla \delta \vec u^o + (\nabla \delta \vec u^o)^\top ) + \lambda \dv \delta \vec u^o \mathbb I_2 \bigr\rbrack \cdot \delta \vec u_\mrm m \idx \\
	&  ~~ + 2 \int \dv \biggl\lbrack p_{\mrm{m},1} \delta \biggl( \dfrac{\nabla  \vec u^o + (\nabla  \vec u^o)^\top}{\sqrt{|\nabla  \vec u^o + (\nabla  \vec u^o)^\top|^2 + \varepsilon^2} } \biggr)  + p_{\mrm m, 1} \delta \biggl( \dfrac{\dv  \vec u^o \mathbb I_2 }{\sqrt{|\dv \vec u^o|^2 + \varepsilon^2}}\biggr) \biggr\rbrack \cdot \delta \vec u_\mrm m \idx \\
	& ~~ + 2 \int \delta p_\mrm m \biggl\lbrack \dfrac{\nabla  \vec u^o_2 + (\nabla  \vec u^o_2)^\top}{\sqrt{|\nabla  \vec u^o_2 + (\nabla  \vec u^o_2)^\top|^2+ \varepsilon^2 }} + \dfrac{\dv  \vec u^o_2 \mathbb I_2 }{\sqrt{|\dv \vec u^o_2|^2 + \varepsilon^2}} \biggr\rbrack : \nabla \delta \vec u_\mrm m \idx.
\end{align*}
Therefore, applying H\"older's inequality and the Sobolev embedding inequality implies
\begin{align*}
	& \mathcal R_{22} \lesssim C_{\varepsilon,\mu,\lambda} (1 + \overline h + \norm{\nabla h_{\mrm m,1},\nabla A_{\mrm m, 1}}{\Lnorm{4}} )\norm{\delta\vec u^o}{\Hnorm{2}} \norm{\delta \vec u_m}{\Lnorm{2}} \\
	& ~~~~ + C_\varepsilon \overline h \sum_{j=1}^2\norm{\nabla^2 \vec u^o_j}{\Lnorm{2}} \norm{\nabla \delta \vec u^o}{\Lnorm{4}} \norm{\delta \vec u_m}{\Lnorm{2}}^{1/2}\norm{\delta \vec u_m}{\Hnorm{1}}^{1/2}\\
	& ~~~~ + (1 + \overline h)\norm{\delta h_\mrm m,\delta A_\mrm m}{\Lnorm{2}} \norm{\nabla \delta \vec u_\mrm m}{\Lnorm{2}},
\end{align*}
where we have used the identity
\begin{align*}
	& \delta\biggl( \dfrac{g}{\sqrt{|g|^2+ \varepsilon^2}} \biggr) = \dfrac{\delta g}{\sqrt{|g_1|^2+ \varepsilon^2} } \\
	& \quad -\dfrac{g_2 \delta \abs{g}{2}}{\sqrt{|g_1|^2+\varepsilon^2}\sqrt{|g_2|^2+\varepsilon^2}(\sqrt{|g_1|^2+\varepsilon^2}+\sqrt{|g_2|^2+\varepsilon^2})}
\end{align*}
for  $ g = \nabla \vec u^o + (\nabla \vec u^o)^\top $ and $ \dv \vec u^o\mathbb I_2 $, respectively.

Then, after substituting the bounds in \eqref{est:h-A-lin-total} and \eqref{est:u-lin-bound} and applying interpolation inequalities, one can obtain from \eqref{est:u-diff-lin-001} that
\begin{equation}\label{est:u-diff-lin-002}
	\begin{aligned}
		& \rho_\mrm{ice} \dfrac{d}{dt} \norm{h_{\mrm m,1}^{1/2} \delta \vec u_\mrm m}{\Lnorm{2}}^2 +  \iota \norm{\delta \vec u_\mrm m}{\Hnorm{2}}^2 \leq C_{\sigma,\varepsilon,\mu,\lambda, \mathfrak c_o, \mathfrak c_\mrm{in}}\norm{\delta \vec u_\mrm m}{\Lnorm{2}}^2 \\
		& ~~~~ + C_{\overline h} (1 + \norm{\dt \vec u_{\mrm m,2}}{\Lnorm{2}}^2) ( \norm{\delta h_\mrm m}{\Lnorm{4}}^2 + \norm{\delta h_\mrm m, \delta A_\mrm m}{\Lnorm{2}}^2 ) \\
		& ~~~~ + \sigma \norm{\delta \vec u^o}{\Hnorm{2}}^2,
	\end{aligned}
\end{equation}
where Young's inequality is applied.

Thus, after substituting \eqref{est:h-A-diff-lin-002} into \eqref{est:u-diff-lin-002} and applying Gr\"onwall's inequality to the resultant, one has
\begin{align*}
	& \sup_{0\leq s\leq t} \norm{\delta \vec u_\mrm m (s)}{\Lnorm{2}}^2 + \int_0^t \norm{\delta \vec u_\mrm m (s)}{\Hnorm{2}}^2 \,ds \leq \sigma C_{\iota,\underline h, \overline h, \varepsilon,\omega,\nu,\mathfrak c_o, \mathfrak c_{\mrm{in}}}\\
	& ~~~~ \times \exp \biggl\lbrack C_{\sigma,\iota,\mu,\lambda,\underline h, \overline h, \varepsilon,\omega,\nu,\mathfrak c_o, \mathfrak c_{\mrm{in}}} ( t  + t^{2}) \biggr\rbrack   \int_0^t \norm{\delta \vec u^o}{\Hnorm{2}}^2 \,ds.
\end{align*}
Therefore, after choosing $ \sigma $ and $ t $ small enough, one can conclude that
\begin{equation}\label{est:contracting}
	\begin{aligned}
	& \sup_{0\leq s\leq t} \norm{\delta \vec u_\mrm m (s)}{\Lnorm{2}}^2 + \int_0^t \norm{\delta \vec u_\mrm m (s)}{\Hnorm{2}}^2 \,ds \\
	& ~~~ \leq \dfrac{1}{2} \biggl( \sup_{0\leq s\leq t} \norm{\delta \vec u^o (s)}{\Lnorm{2}}^2 + \int_0^t \norm{\delta \vec u^o (s)}{\Hnorm{2}}^2 \,ds \biggr).
	\end{aligned}
\end{equation}
Now we update the smallness of $ T^* $, so that \eqref{est:contracting} holds true for $ t \in (0,T^*] $. Then the map $ \mathfrak M $, defined in \eqref{def:mapping}, is contracting with constant $ 1/2 $. By means of Banach's fixed point theorem, we conclude that there exists a unique solution to \eqref{eq:sea-ice-app} in $ \mathfrak X_\mrm {T^*} $.

What is left is to show that such solutions are stable. Namely, they continuously depend on the initial data. Let $  (\vec u_j, h_j, A_j ) $ be two solutions to \eqref{eq:sea-ice-app}, associated with initial data $  (\vec u_{\mrm{in}, j}, h_{\mrm{in}, j}, A_{\mrm{in}, j} ) $, $ j = 1,2 $, satisfying \eqref{def:initial-bound}. Then it is easy to check that \eqref{est:h-A-diff-lin-002} and \eqref{est:contracting} still hold true with $\delta \vec u^o, \delta \vec u_m, \delta h_\mrm m, \delta A_\mrm m $ replaced by $ \delta \vec u := \vec u_1 - \vec u_2,  \delta h:= h_1 - h_2, \delta A:= A_1 - A_2 $, with additional initial data on the righthand side, i.e.,
\begin{equation}\label{est:stability-app}
	\begin{aligned}
		& \sup_{0\leq s\leq t} \bigl( \norm{\delta h(s), \delta A(s)}{\Lnorm{4}}^2 + \norm{\delta \vec u(s)}{\Lnorm{2}}^2 \bigr) + \int_0^t \norm{\delta \vec u(s)}{\Hnorm{2}}^2 \,ds \\
		& ~~ \leq C_{\varepsilon,\omega,\nu,\underline h,\overline h, \mathfrak c_o,\mathfrak c_\mrm{in}} \bigl( \norm{h_{\mrm{in},1} - h_{\mrm{in}, 2}, A_{\mrm{in},1} - A_{\mrm{in}, 2}}{\Lnorm{4}}^2 + \norm{\vec u_{\mrm{in},1} - \vec u_{\mrm{in}, 2}}{\Lnorm{2}}^2 \bigr).
	\end{aligned}
\end{equation}

Hence, we have established the local-in-time well-posedness of strong solutions to system \eqref{eq:sea-ice-app}. We would like to remind readers that the estimates obtained in this section depend on $ (\mu, \lambda,\iota, \nu) $. In the next section, we aim at removing such dependency.


\section{Well-posedness of solutions to \eqref{eq:sea-ice-app-2} with
$ \underline{h}> 0 $}
\label{sec:well-posedness}

\subsection{$ (\mu, \lambda,\iota, \nu) $-independent estimates of solutions to \eqref{eq:sea-ice-app}}

We shall only present the uniform-in-$ (\mu, \lambda,\iota,\nu) $ {\it a priori} estimate in this subsection, based on which the standard different quotient argument can be established.

Throughout this section, we use the notation $ X \lesssim Y $ to represent $ X \leq C Y $ for some generic constant $ C \in (0,\infty) $, which may be different from line to line, and depend on $ \varepsilon,  \omega, \underline h, \overline h $, but is independent of $ (\mu, \lambda,\iota,\nu) $.

To begin with, let
\begin{equation}\label{def:uni-energy-total}
	\mathcal E(t) := \sup_{0\leq s \leq t} \norm{\vec u(s),h(s),A(s)}{\Hnorm{3}}^2 + \int_0^t\norm{\vec u(s)}{\Hnorm{4}}^2 \,ds,
\end{equation}
and
\begin{equation}\label{def:uni-energy-total-2}
    \begin{aligned}
    & \mathfrak E(t) := \sup_{0\leq s \leq t} \norm{\vec u(s),h(s),A(s)}{\Hnorm{3}}^2 \\
    &\quad + \int_0^t \int \biggl( \dfrac{|\nabla^3(\nabla\vec u(s) + \nabla \vec u^\top(s))|^2}{(|\nabla\vec u(s) + \nabla \vec u^\top(s)|^2 +\varepsilon^2)^{3/2}}+ \dfrac{|\nabla^3\dv \vec u(s)|^2}{(|\dv \vec u(s)|^2 + \varepsilon^2)^{3/2}} \biggr) \idx \,ds.
    \end{aligned}
\end{equation}
One can easily check that $ \mathcal E $ and $ \mathfrak E $ are essentially equivalent in the sense that estimates on one imply estimates on the other. Indeed, it is trivial that $ \mathfrak E \lesssim \mathcal E $. On the other hand, applying integration by parts yields that
\begin{equation}\label{eqvl:uni-h4}
\begin{aligned}
    & \int |\nabla^4 \vec u|^2 \idx = \dfrac{1}{2} \int |\nabla^3(\nabla \vec u+ \nabla \vec u^\top)|^2 \idx - \int |\nabla^3 \dv \vec u|^2 \idx \\
    &\quad \lesssim (\varepsilon^3 + \norm{\vec u}{\Hnorm{3}}^3) \int \biggl( \dfrac{|\nabla^3(\nabla\vec u + \nabla \vec u^\top)|^2}{(|\nabla\vec u + \nabla \vec u^\top|^2 +\varepsilon^2)^{3/2}}+ \dfrac{|\nabla^3\dv \vec u|^2}{(|\dv \vec u|^2 + \varepsilon^2)^{3/2}} \biggr) \idx.
\end{aligned}
\end{equation}
Therefore, we have
\begin{equation}\label{eqvl:uni-energy}
    \mathfrak E(t) \lesssim \mathcal E(t) \lesssim (1 + t + \mathfrak E^2(t))\mathfrak E(t).
\end{equation}

\subsubsection*{Estimates for $ h $ and $ A $}
It is easy to check that \eqref{eq:positivity-A-app}, \eqref{eq:upper-A-app}, \eqref{eq:positivity-h-app}, \eqref{eq:upper-h-app}, and \eqref{eq:lower-h-app} also hold true with $ A_\mrm m, h_\mrm m, \vec u^o $ replaced by $ A, h, \vec u $, respectively. Therefore, for $ s \in(0,t) $ with $ t $ satisfying \eqref{cnd:small-time-001}, with $ \vec u^o $ replaced by $ \vec u $, we have
\begin{equation}\label{unest:ptwbd-h-A}
    0 \leq A \leq 1, \qquad 0 < \dfrac{1}{4} \underline h \leq h \leq 4 \overline h.
\end{equation}
Notice that the smallness of $ t $ here is independent of $ (\mu,\lambda, \iota, \nu ) $.

Next, we shall establish the regularity estimates of $ A $ and $ h $. Indeed, after applying $ \partial^3 $ to \eqref{eq:sea-ice-app=2} and \eqref{eq:sea-ice-app=3}, one can obtain the following equations:
\begin{subequations}
\begin{gather}
\begin{gathered}
\dt\partial^3 h + \vec u \cdot \nabla \partial^3 h = \partial^3 \mathcal S_{h,\mu,\nu} - \partial^3 ( h \dv \vec u) \\
+ \bigl( \vec u \cdot \nabla \partial^3 h - \partial^3( \vec u \cdot \nabla h ) \bigr), \end{gathered} \label{eq:d-3-h} \\
\begin{gathered}
\dt\partial^3 A + \vec u \cdot \nabla \partial^3 A = \partial^3 \mathcal S_{A,\omega,\nu} + \partial^3 (A \dv \vec u \cdot \chi^\omega_A) \\
- \partial^3( A \dv \vec u) + \bigl( \vec u \cdot \nabla \partial^3 A - \partial^3 (\vec u \cdot \nabla A) \bigr).\end{gathered} \label{eq:d-3-A}
\end{gather}
\end{subequations}

Then, applying the $ L^2 $-inner product of \eqref{eq:d-3-h} and \eqref{eq:d-3-A} with $ 2 \partial^3 h $ and $ \partial^3 A $, respectively, and integration by parts in the resultant leads to
\begin{subequations}
\begin{gather}
    \label{unest:h-001} \begin{aligned}
    & \dfrac{d}{dt} \norm{\partial^3 h}{\Lnorm{2}}^2 = \underbrace{\int \bigl( \dv \vec u |\partial^3 h|^2 - 2 \partial^3 (h\dv \vec u) \partial^3 h \bigr) \idx}_{\mathcal I_{1}}\\
    & \quad + \underbrace{2 \int \bigl( \vec u \cdot \nabla \partial^3 h - \partial^3( \vec u \cdot \nabla h ) \bigr) \partial^3 h \idx}_{\mathcal I_{2}} + \underbrace{2 \int \partial^3 \mathcal S_{h,\mu,\nu} \partial^3 h \idx}_{\mathcal I_{3}},
    \end{aligned}\\
    \label{unest:A-001} \begin{aligned}
    & \dfrac{d}{dt} \norm{\partial^3 A}{\Lnorm{2}}^2 = \underbrace{\int \bigl( \dv \vec u |\partial^3 A|^2 - 2 \partial^3(A \dv \vec u) \partial^3 A \bigr) \idx}_{\mathcal I_{4}} \\
    & \quad + \underbrace{2 \int \bigl( \vec u \cdot \nabla \partial^3 A - \partial^3 (\vec u \cdot \nabla A) \bigr) \idx}_{\mathcal I_{5}} + \underbrace{2 \int \partial^3 \mathcal S_{A,\mu,\nu} \partial^3 A \idx}_{\mathcal I_{6}} \\
    & \quad + \underbrace{2 \int \partial^3 ( A \dv \vec u \cdot \chi^\omega_A ) \partial^3 A \idx}_{\mathcal I_{7}}.
    \end{aligned}
\end{gather}
\end{subequations}

Directly applying H\"older's inequality and the Sobolev embedding inequality leads to the following estimates:
\begin{equation}\label{unest:h-A-001}
\begin{gathered}
     \mathcal I_1 + \mathcal I_2 + \mathcal I_4 + \mathcal I_5 + \mathcal I_7 \lesssim \mathcal H( \norm{\vec u, h, A}{\Hnorm{3}} ) \\
    + \norm{\vec u}{\Hnorm{4}} \norm{h,A}{\Hnorm{3}}^2 .
\end{gathered}
\end{equation}
Similarly,
\begin{equation}\label{unest:h-A-002}
    \mathcal I_3 + \mathcal I_6 \lesssim \mathcal H(\norm{h,A}{\Hnorm{3}}).
\end{equation}

Therefore, after substituting estimates \eqref{unest:h-A-001} and \eqref{unest:h-A-002} into \eqref{unest:h-001} and \eqref{unest:A-001}, one can derive that
\begin{gather*}
    \dfrac{d}{dt} \norm{\partial^3 h,\partial^3 A}{\Lnorm{2}}^2  \lesssim \mathcal H( \norm{\vec u,h,A}{\Hnorm{3}} )
    + \norm{\vec u}{\Hnorm{4}}  \norm{h,A}{\Hnorm{3}}^2 .
\end{gather*}
Similar estimates also hold for lower order derivatives. Hence we have shown that
\begin{gather*}
    \dfrac{d}{dt} \norm{h,A}{\Hnorm{3}}^2 \leq \mathcal H( \norm{\vec u,h,A}{\Hnorm{3}} )
    + C_{\omega,\underline h, \overline h} \norm{\vec u}{\Hnorm{4}} \norm{h,A}{\Hnorm{3}}^2 ,
\end{gather*}
for some constant $ C_{\omega,\underline h, \overline h} \in (0,\infty) $, independent of $ \iota $ and $ \nu $. Consequently, applying Gr\"onwall's inequality concludes that
\begin{equation}\label{unest:h-A-101}
    \begin{gathered}
        \sup_{0\leq s \leq t}  \norm{h(s),A(s)}{\Hnorm{3}}^2  \leq e^{C_{\omega,\underline h, \overline h} \int_0^t \norm{\vec u(s)}{\Hnorm{4}}\,ds} \\
        \times \biggl( \norm{h_\mrm{in},A_\mrm{in}}{\Hnorm{3}}^2 + \int_0^t \mathcal H(\norm{\vec u(s),h(s),A(s)}{\Hnorm{3}} )\,ds\biggr).
    \end{gathered}
\end{equation}

\subsubsection*{Estimates for $ \vec u $}

After applying $ \partial^3 $ to  \eqref{eq:sea-ice-app=1}, one can obtain the following equation:
\begin{equation}\label{eq:d-3-u}
\begin{gathered}
    m ( \dt \partial^3 \vec u + \vec u \cdot \nabla \partial^3 \vec u ) + \nabla \partial^3 p = \dv \partial^3 \mathbb S_{\varepsilon} + \dv \partial^3 \mathbb S_{\mu, \lambda} \\
    - \iota \Delta^2 \partial^3 \vec u
    + \partial^3 \mathcal F
    + \lbrack m  \dt \partial^3 \vec u - \partial^3 (m \dt \vec u) \rbrack \\
    + \lbrack m \vec u \cdot \nabla \partial^3 \vec u - \partial^3 ( m \vec u \cdot \nabla \vec u) \rbrack.
\end{gathered}
\end{equation}

Then, applying the $ L^2 $-inner product of \eqref{eq:d-3-u} with $ 2 \partial^3 \vec u $ and integration by parts in the resultant leads to
\begin{equation}\label{unest:u-001}
    \begin{aligned}
    &  \dfrac{d}{dt} \norm{\rho_\mrm{ice}^{1/2} h^{1/2} \partial^3 \vec u}{\Lnorm{2}}^2 + 2 \mu \norm{\nabla \partial^3 \vec u}{\Lnorm{2}}^2 + 2(\mu +\lambda) \norm{\dv \partial^3 \vec u}{\Lnorm{2}}^2 \\
    & \qquad + 2\iota \norm{\nabla^2 \partial^3 \vec u}{\Lnorm{2}}^2 = \underbrace{\int \lbrack \rho_\mrm{ice} \dt h + \dv ( \rho_\mrm{ice} h \vec u) \rbrack \abs{\partial^3 u}{2} \idx}_{\mathcal I_8} \\
    & \qquad \underbrace{- 2 \int \partial^3 \mathbb S_\varepsilon : \nabla \partial^3 \vec u \idx}_{\mathcal I_9} + \underbrace{2 \int \lbrack m  \dt \partial^3 \vec u - \partial^3 (m \dt \vec u) \rbrack \cdot \partial^3 \vec u \idx}_{\mathcal I_{10}}  \\
    & \qquad + \underbrace{2\int \lbrack \rho_\mrm{ice} h \vec u \cdot \nabla \partial^3 \vec u - \partial^3 ( \rho_\mrm{ice} h \vec u \cdot \nabla \vec u) \rbrack \cdot \partial^3 \vec u \idx}_{\mathcal I_{11}} \\
    & \qquad + \underbrace{2 \int \partial^3 p \dv \partial^3 \vec u \idx}_{\mathcal I_{12}} \underbrace{- 2 \int \partial^2 \mathcal F \cdot \partial^4 \vec u \idx}_{\mathcal I_{13}}. 
    \end{aligned}
\end{equation}
The estimates of $ \mathcal I_\mrm j $, $ \mrm j \in \lbrace 8,11,12 \rbrace $, are standard, which we will record below.
Applying H\"older's inequality and the Sobolev embedding inequality yields  that
\begin{equation}\label{unest:u-002}
\begin{aligned}
    \mathcal I_8 & \lesssim \bigl( \norm{\dt h}{\Lnorm{2}} + \norm{\dv (h \vec u)}{\Lnorm{2}} \bigr)  \norm{\partial^3 \vec u}{\Lnorm{2}} \norm{\partial^3 \vec u}{\Hnorm{1}} \\
    & \lesssim \bigl( \norm{h}{\Lnorm{\infty}} + \norm{\nabla h}{\Lnorm{4}} \bigr)  \norm{\vec u}{\Hnorm{3}}^2 \norm{\vec u}{\Hnorm{4}},\\
    \mathcal I_{11} & \lesssim \norm{h}{\Hnorm{3}} \norm{\vec u}{\Hnorm{3}}^2 \norm{\vec u}{\Hnorm{4}}, \\
    \mathcal I_{12} & \lesssim  \bigl( \norm{A}{\Hnorm{3}}^3 + 1 \bigr)\norm{h}{\Hnorm{3}}  \norm{\vec u}{\Hnorm{4}}.
\end{aligned}
\end{equation}

To estimate $ \mathcal I_{13} $, notice that
$$
    \norm{\partial^2 \mathcal F}{\Lnorm{2}} \lesssim \norm{\partial^2 \bigl( |\vec U_\mrm w - \vec u|(\vec U_\mrm w - \vec u) \bigr)}{\Lnorm{2}} + \norm{h}{\Hnorm{2}}\norm{\vec u}{\Hnorm{2}} + \mrm{l.o.t},
$$
where $ \mrm{l.o.t} $ represents lower order terms of $ \vec u $. Direct calculation yields that
\begin{gather*}
    \partial^2 \bigl( |\vec U_\mrm w - \vec u|(\vec U_\mrm w - \vec u) \bigr) =  |\vec U_\mrm w - \vec u|\partial^2 (\vec U_\mrm w - \vec u) \\
    + 2 \dfrac{(\vec U_\mrm w - \vec u) \cdot \partial (\vec U_\mrm w - \vec u)}{|\vec U_\mrm w - \vec u|} \partial(\vec U_\mrm w - \vec u) \\
    + \biggl( \dfrac{(\vec U_\mrm w - \vec u) \cdot \partial^2 (\vec U_\mrm w - \vec u) + |\partial (\vec U_\mrm w - \vec u)|^2}{|\vec U_\mrm w - \vec u|} \\
    - \dfrac{\bigl((\vec U_\mrm w - \vec u) \cdot \partial (\vec U_\mrm w - \vec u)\bigr)^2}{|\vec U_\mrm w - \vec u|^3} \biggr) (\vec U_\mrm w - \vec u),
\end{gather*}
which implies
$$
\norm{\partial^2 \bigl( |\vec U_\mrm w - \vec u|(\vec U_\mrm w - \vec u) \bigr)}{\Lnorm{2}} \lesssim \norm{\vec U_\mrm w - \vec u}{\Hnorm{2}}^2 + \norm{\vec U_\mrm w - \vec u}{\Hnorm{2}}^3.
$$
Therefore, we have
\begin{equation}\label{unest:u-003}
    \mathcal I_{13} \lesssim \norm{\partial^2 \mathcal F}{\Lnorm{2}} \norm{\partial^4\vec u}{\Lnorm{2}} \lesssim \bigl(\norm{\vec u}{\Hnorm{3}}^3 + \norm{h}{\Hnorm{2}}^2 + 1  \bigr) \norm{u}{\Hnorm{4}}.
\end{equation}

In order to estimate $ \mathcal I_{10} $, we first rewrite $ \mathcal I_{10} $ as follows,
\begin{equation}\label{unest:u-0031}
    \mathcal I_{10} = - 2 \int \partial^3 m \dt \vec u \cdot \partial^3 \vec u \idx
    + 6 \int \partial m \dt \partial \vec u \cdot \partial^4 \vec u \idx,
\end{equation}
where we have applied integration by parts. Next, we will use equation \eqref{eq:sea-ice-app=1} to substitute $ \dt \vec u $ and $ \dt \partial \vec u $ in \eqref{unest:u-0031}.
Indeed, after rearranging \eqref{eq:sea-ice-app=1}, it follows
\begin{gather*}
    \dt \vec u = \dfrac{\dv \mathbb S_{\varepsilon,\mu,\lambda}}{m}  + \dfrac{\mathcal F}{m} - \dfrac{\nabla p}{m} - \vec u \cdot\nabla \vec u - \iota\dfrac{\Delta^2 \vec u}{m},\\
    \begin{aligned}
    \dt \partial \vec u = & \dfrac{\dv \partial \mathbb S_{\varepsilon,\mu,\lambda}}{m} - \dfrac{\dv \mathbb S_{\varepsilon,\mu,\lambda}}{m^2} \partial m
    + \dfrac{\partial \mathcal F}{m} - \dfrac{\mathcal F}{m^2} \partial m \\
    & - \dfrac{\nabla\partial p}{m} + \dfrac{\nabla p}{m^2} \partial m
     - \partial \vec u\cdot \nabla \vec u - \vec u \cdot \nabla \partial \vec u \\
     & - \iota \dfrac{\Delta^2 \partial \vec u}{m} + \iota \dfrac{\Delta^2 \vec u}{m^2} \partial m.
     \end{aligned}
\end{gather*}
Then similarly as before, directly applying H\"older's inequality and the Sobolev embedding inequality leads to,
\begin{gather*}
\norm{\dt \vec u}{\Lnorm{4}} +
\norm{\dt \partial \vec u}{\Lnorm{2}} \lesssim \mathcal H (\norm{\vec u}{\Hnorm{3}}, \norm{A}{\Hnorm{2}}, \norm{h}{\Hnorm{2}}) \\
+ \iota (1 + \norm{h}{\Hnorm{2}}) \norm{\vec u}{\Hnorm{5}}.
\end{gather*}
Therefore, one can derive that,
\begin{equation}\label{unest:u-004}
\begin{aligned}
    \mathcal I_{10} \lesssim & \norm{\partial^3 m}{\Lnorm{2}} \norm{\dt \vec u}{\Lnorm{4}} \norm{\partial^3 \vec u}{\Lnorm{4}} \\ & \qquad + \norm{\partial m}{\Lnorm{\infty}} \norm{\dt \partial \vec u}{\Lnorm{2}} \norm{\partial^4 \vec u}{\Lnorm{2}} \\
    \lesssim & \mathcal H (\norm{\vec u}{\Hnorm{3}}, \norm{A}{\Hnorm{2}}, \norm{h}{\Hnorm{3}}) \norm{\vec u}{\Hnorm{4}} \\
    & \qquad + \iota (\norm{h}{\Hnorm{3}} + \norm{h}{\Hnorm{3}}^2) \norm{\vec u}{\Hnorm{5}} \norm{\vec u}{\Hnorm{4}}.
\end{aligned}
\end{equation}

Lastly, we will estimate $ \mathcal I_9 $. Notice that,
\begin{gather*}
    \mathcal I_9 = 
    - \int \partial^3 \biggl( p \dfrac{\nabla \vec u + \nabla \vec u^\top}{\sqrt{|\nabla \vec u+ \nabla \vec u^\top|^2 + \varepsilon^2}} \biggr) : \partial^3 (\nabla \vec u + \nabla  \vec u^\top ) \idx \\
    - 2 \int \partial^3 \biggl(p \dfrac{\dv \vec u}{\sqrt{|\dv \vec u|^2 + \varepsilon^2}}\biggr) \partial^3 \dv \vec u \idx.
\end{gather*}
Denote by $ \mrm{D} \vec u \in \lbrace \nabla \vec u + \nabla \vec u^\top, \dv \vec u  \rbrace  $.
In this notation, estimating $ \mathcal I_9 $ amounts to determining an estimate for
$$
    \int \partial^3 \biggl(p \dfrac{\mrm D \vec u}{\sqrt{|\mrm D \vec u|^2 + \varepsilon^2}}\biggr) \cdot \partial^3 \mrm D \vec u \idx.
$$
Direct calculation shows that
\begin{gather*}
    \int \partial^3 \biggl(p \dfrac{\mrm D \vec u}{\sqrt{|\mrm D \vec u|^2 + \varepsilon^2}}\biggr) \cdot \partial^3 \mrm D \vec u \idx = \int p \biggl( \dfrac{|\partial^3 \mrm D \vec u|^2}{\sqrt{|\mrm D \vec u|^2 + \varepsilon^2}}  -  \dfrac{(\mrm D \vec u \cdot \partial^3 \mrm D \vec u)^2}{(|\mrm D \vec u|^2 + \varepsilon^2)^{3/2}} \biggr) \idx \\
    - \underbrace{3 \int p \dfrac{(\mrm D \vec u \cdot \partial \mrm D \vec u)(\partial^2 \mrm D \vec u \cdot \partial^3 \mrm D \vec u)+(\mrm D \vec u \cdot \partial^2 \mrm D \vec u)(\partial \mrm D \vec u \cdot \partial^3 \mrm D \vec u)}{(|\mrm D \vec u|^{2} + \varepsilon^2)^{3/2}}\idx}_{\mathcal L_1} \\
    -\underbrace{3 \int p \dfrac{(\partial \mrm D \vec u \cdot \partial^2 \mrm D \vec u)(\mrm D \vec u \cdot \partial^3 \mrm D \vec u)}{(|\mrm D \vec u|^2 + \varepsilon^2)^{3/2}} \idx}_{\mathcal L_2} \\
    + \underbrace{9 \int p \dfrac{(\mrm D \vec u \cdot \partial \mrm D \vec u) ( \mrm D \vec u \cdot \partial^2 \mrm D \vec u) ( \mrm D \vec u \cdot \partial^3 \mrm D \vec u)}{(|\mrm D \vec u|^2 + \varepsilon^2)^{5/2}} \idx}_{\mathcal L_3}\\
    - \underbrace{3 \int p \dfrac{|\partial \mrm D \vec u|^2 ( \partial \mrm D \vec u \cdot \partial^3 \mrm D \vec u)}{(|\mrm D \vec u|^2 + \varepsilon^2)^{3/2}} \idx}_{\mathcal L_4}
    \\
    + \underbrace{9 \int p \dfrac{(\mrm D \vec u \cdot \partial \mrm D \vec u)^2 ( \partial \mrm D \vec u \cdot \partial^3 \mrm D \vec u) + |\partial \mrm D \vec u|^2 ( \mrm D \vec u \cdot \partial \mrm D \vec u) (\mrm D \vec u \cdot \partial^3 \mrm D \vec u)}{(|\mrm D \vec u|^2 + \varepsilon^2)^{5/2}}\idx}_{\mathcal L_5}\\
    - \underbrace{15 \int p \dfrac{(\mrm D \vec u\cdot \partial \mrm D \vec u)^3 (\mrm D \vec u \cdot \partial^3 \mrm D \vec u)}{(|\mrm D \vec u|^2 + \varepsilon^2)^{7/2}} \idx}_{\mathcal L_6} \\
    + \underbrace{3 \int \biggl\lbrack \partial p \partial^2 \biggl( \dfrac{\mrm D \vec u}{\sqrt{|\mrm D\vec u|^2 + \varepsilon^2}} \biggr) \cdot \partial^3 \mrm D \vec u + \partial^2 p \partial \biggl( \dfrac{\mrm D \vec u}{\sqrt{|\mrm D\vec u|^2 + \varepsilon^2}} \biggr) \cdot \partial^3 \mrm D \vec u \biggr\rbrack \idx}_{\mathcal L_7}\\
    + \underbrace{\int \partial^3 p  \dfrac{\mrm D \vec u \cdot \partial^3 \mrm D \vec u}{\sqrt{|\mrm D\vec u|^2 + \varepsilon^2}}  \idx}_{\mathcal L_8}.
\end{gather*}
Notice that
\begin{gather*}
    \dfrac{|\partial^3 \mrm D \vec u|^2}{\sqrt{|\mrm D \vec u|^2 + \varepsilon^2}}  -  \dfrac{(\mrm D \vec u \cdot \partial^3 \mrm D \vec u)^2}{(|\mrm D \vec u|^2 + \varepsilon^2)^{3/2}} \geq 
    \varepsilon^2 \dfrac{|\partial^3 \mrm D \vec u|^2}{(|\mrm D \vec u|^2 + \varepsilon^2)^{3/2}} .
\end{gather*}
Therefore, applying H\"older's inequality and the Sobolev embedding inequality implies that,
\begin{gather*}
	|\mathcal L_4|+|\mathcal L_5|+|\mathcal L_6|+|\mathcal L_7|+|\mathcal L_8| \lesssim \norm{p}{\Hnorm{3}} (1 + \norm{\mrm D \vec u}{\Hnorm{2}}^3) \norm{\partial^3 \mrm D \vec u}{\Lnorm{2}},  \\
    |\mathcal L_1|+|\mathcal L_2|+|\mathcal L_3| \lesssim \norm{p}{\Lnorm{\infty}}\norm{\mrm D\vec u}{\Hnorm{2}}^{3/2} \norm{\mrm D \vec u}{\Hnorm{3}}^{3/2}.
\end{gather*}
Therefore,
\begin{equation*}
\begin{gathered}
    \int \partial^3 \biggl(p \dfrac{\mrm D \vec u}{\sqrt{|\mrm D \vec u|^2 + \varepsilon^2}}\biggr) \cdot \partial^3 \mrm D \vec u \idx \geq \varepsilon^2 \int \dfrac{p |\partial^3 \mrm D \vec u|^2}{(|\mrm D \vec u|^2 + \varepsilon^2)^{3/2}} \idx  \\
    - \norm{p}{\Hnorm{3}} (1 + \norm{\mrm D \vec u}{\Hnorm{2}}^3) \norm{\partial^3 \mrm D \vec u}{\Lnorm{2}} \\
    - \norm{p}{\Lnorm{\infty}}\norm{\mrm D\vec u}{\Hnorm{2}}^{3/2} \norm{\mrm D \vec u}{\Hnorm{3}}^{3/2}.
\end{gathered}
\end{equation*}
Thus, we have shown that, thanks to the fact $ p \geq  c_p \underline h/4 > 0 $,
\begin{equation}\label{unest:u-005}
    \begin{gathered}
    \mathcal I_9 \leq - \dfrac{\varepsilon^2 c_p \underline h }{4} \int \biggl(\dfrac{|\partial^3 (\nabla \vec u + \nabla \vec u^\top)|^2}{(|(\nabla \vec u + \nabla \vec u^\top)|^2 + \varepsilon^2)^{3/2}} + 2\dfrac{|\partial^3 \dv \vec u|^2}{(|\dv \vec u|^2 + \varepsilon^2)^{3/2}}\biggr) \idx \\
    + \mathcal H (\norm{\vec u,A,h}{\Hnorm{3}}) \norm{\vec u}{\Hnorm{4}} + \norm{\vec u}{\Hnorm{3}}^{3/2}\norm{\vec u}{\Hnorm{4}}^{3/2}.
    \end{gathered}
\end{equation}
In addition, notice that, according to \eqref{eqvl:uni-h4},
\begin{equation}\label{unest:u-006}
    \begin{aligned}
    & \norm{\vec u}{\Hnorm{4}} \lesssim  \norm{\vec u}{\Hnorm{3}} + \norm{\nabla^4 \vec u}{\Lnorm{2}}
    \lesssim \norm{\vec u}{\Hnorm{3}} \\
    & \quad + \biggl\lbrack(\varepsilon^3 + \norm{\vec u}{\Hnorm{3}}^3) \int \biggl( \dfrac{|\nabla^3(\nabla\vec u + \nabla \vec u^\top)|^2}{(|\nabla\vec u + \nabla \vec u^\top|^2 +\varepsilon^2)^{3/2}}+ \dfrac{|\nabla^3\dv \vec u|^2}{(|\dv \vec u|^2 + \varepsilon^2)^{3/2}} \biggr) \idx\biggr\rbrack^{1/2}.
    \end{aligned}
\end{equation}

To sum up, after substituting estimates \eqref{unest:u-002}, \eqref{unest:u-003}, \eqref{unest:u-004}, \eqref{unest:u-005},and \eqref{unest:u-006} into \eqref{unest:u-001}, and applying Young's inequality, one can derive that
\begin{gather*}
    \dfrac{d}{dt} \norm{\rho_\mrm{in}^{1/2} h^{1/2} \partial^3 \vec u}{\Lnorm{2}}^2
    + 2\iota \norm{\nabla^2 \partial^3 \vec u}{\Lnorm{2}}^2 - 2 \iota^2 \norm{\nabla^5 \vec u}{\Lnorm{2}}^2 \\
    + \dfrac{\varepsilon^2 c_p \underline h }{8} \int \biggl( \dfrac{|\partial^3 (\nabla \vec u + \nabla \vec u^\top)|^2}{(|(\nabla \vec u + \nabla \vec u^\top)|^2 + \varepsilon^2)^{3/2}} + 2 \dfrac{|\partial^3 \dv \vec u|^2}{(|\dv \vec u|^2 + \varepsilon^2)^{3/2}}\biggr) \idx \\
    \leq \mathcal H(\norm{\vec u,A,h}{\Hnorm{3}},\iota),
\end{gather*}
which implies, recalling $ \iota \in (0,1) $,
\begin{equation*}
\begin{gathered}
    \sup_{0\leq s\leq t} \norm{\nabla^3 \vec u(s)}{\Lnorm{2}}^2
    + (\iota - \iota^2) \int_0^t\norm{\nabla^5 \vec u(s)}{\Lnorm{2}}^2 \,ds \\ + \int_0^t \int \biggl( \dfrac{|\partial^3 (\nabla \vec u(s) + \nabla \vec u^\top(s))|^2}{(|(\nabla \vec u(s) + \nabla \vec u^\top(s))|^2 + \varepsilon^2)^{3/2}} + 2 \dfrac{|\partial^3 \dv \vec u(s)|^2}{(|\dv \vec u(s)|^2 + \varepsilon^2)^{3/2}}\biggr) \idx \,ds \\
    \leq C_{\varepsilon,\underline h, \overline h} \norm{\nabla^3 \vec u_\mrm{in}}{\Lnorm{2}}^2
    + \int_0^t \mathcal H(\norm{\vec u(s),A(s),h(s)}{\Hnorm{3}},\iota)\,ds,
\end{gathered}
\end{equation*}
for some constant $ C_{\varepsilon,\underline h, \overline h} \in (0,\infty) $, independent of $ \mu $, $ \lambda $, $ \iota $, and $ \nu $.

Similar estimates also hold for lower order derivatives. Thus one can conclude that, for $ \iota \ll 1 $ small enough,
\begin{equation}\label{unest:u-101}
\begin{gathered}
     \sup_{0\leq s\leq t} \norm{\vec u(s)}{\Hnorm{3}}^2 \\
     + \int_0^t \int \biggl( \dfrac{|\partial^3 (\nabla \vec u(s) + \nabla \vec u^\top(s))|^2}{(|(\nabla \vec u(s) + \nabla \vec u^\top(s))|^2 + \varepsilon^2)^{3/2}} 
     + 2 \dfrac{|\partial^3 \dv \vec u(s)|^2}{(|\dv \vec u(s)|^2 + \varepsilon^2)^{3/2}}\biggr) \idx \,ds \\
    \leq C_{\varepsilon,\underline h, \overline h} \norm{\vec u_\mrm{in}}{\Hnorm{3}}^2
    + \int_0^t \mathcal H(\norm{\vec u(s),A(s),h(s)}{\Hnorm{3}})\,ds.
\end{gathered}
\end{equation}

\subsubsection*{Uniform estimates}

The summation of
\eqref{unest:h-A-101} and \eqref{unest:u-101} leads to
\begin{align*}
    & \mathfrak E(t) \leq \biggl( e^{C_{\omega,\underline h, \overline h} t^{1/2} \mathcal E^{1/2}(t)} + C_{\varepsilon,\underline h, \overline h} \biggr)
    \\
    & \qquad \times
    \biggl( \norm{h_\mrm{in},A_\mrm{in},\vec u_\mrm{in}}{\Hnorm{3}}^2 + t \times \mathcal H(\mathfrak E(t)) \biggr)\\
    & \quad \leq \biggl( e^{C_{\omega,\underline h, \overline h} t^{1/2} \bigl\lbrack t^{3/2} + \mathfrak E(t) + \mathfrak E^3(t)\bigr\rbrack^{1/2}} + C_{\varepsilon,\underline h, \overline h} \biggr)
    \\
    & \qquad \times
    \biggl( \norm{h_\mrm{in},A_\mrm{in},\vec u_\mrm{in}}{\Hnorm{3}}^2 + t \times \mathcal H(\mathfrak E(t)) \biggr),
\end{align*}
where we have applied \eqref{eqvl:uni-h4} and Young's inequality in the second inequality.
Consequently, for $ t $ small enough, independent of $ \mu,\lambda ,\iota, \nu $, one can conclude that
\begin{equation}\label{unest:h-A-u-total}
    \mathfrak E(t) \leq C_{\varepsilon, \omega, \underline h, \overline h} \times \norm{h_\mrm{in},A_\mrm{in},\vec u_\mrm{in}}{\Hnorm{3}}^2,
\end{equation}
and, thanks to \eqref{eqvl:uni-h4},
\begin{equation}\label{unest:h-A-u-total-2}
    \mathcal E(t) \leq \mathfrak C_\mrm{in}^2,
\end{equation}
for some constant $ \mathfrak C_\mrm{in} \in (0,\infty) $, depending only on $ \varepsilon, \omega,\underline h, \overline h $, and $$ \norm{h_\mrm{in},A_\mrm{in},\vec u_\mrm{in}}{\Hnorm{3}}. $$
Thus we have established the $ (\mu,\lambda, \iota,\nu)$-independent estimates. Therefore, together with the well-posedness theory in Section \ref{sec:first-layer-app} and continuity arguments, the existence time of solutions to \eqref{eq:sea-ice-app} can be extended to some $ T^{**} \in (0,\infty) $, independent of $ (\mu,\lambda, \iota,\nu) $, which might be larger than $ T^* $.

\subsection{Limit as $( \mu,\lambda, \iota, \nu) \rightarrow (0^+,0^+,0^+,0^+) $
}

Denote by $ (\vec u_{\mu,\lambda, \iota,\nu}, h_{\mu,\lambda, \iota,\nu}, A_{\mu,\lambda, \iota,\nu}) $, the solution constructed above to system \eqref{eq:sea-ice-app}.
With \eqref{def:uni-energy-total}, \eqref{unest:h-A-u-total-2},  
and by comparison in system \eqref{eq:sea-ice-app}, it is easy to check that we have the following uniform-in-$(\mu,\lambda, \iota, \nu)$ estimates:
\begin{equation}\label{limit:spatial-regularity}
\begin{gathered}
    \norm{\vec  u_{\mu,\lambda, \iota,\nu}, h_{\mu,\lambda, \iota,\nu}, A_{\mu,\lambda, \iota,\nu}}{L^\infty(0,T^{**};H^3(\Omega))} + \norm{\vec u_{\mu,\lambda, \iota,\nu}}{L^2(0,T^{**};H^4(\Omega))} \\
    + \norm{\dt \vec u_{\mu,\lambda, \iota,\nu},\dt h_{\mu,\lambda, \iota,\nu}, \dt A_{\mu,\lambda, \iota,\nu}}{L^\infty(0,T^{**};L^2(\Omega))} \leq \mathfrak C_\mrm{in},
\end{gathered}
\end{equation}
for some constant $ \mathfrak C_\mrm{in} \in (0,\infty) $, and $ T^{**} \in (0,\infty) $, independent of $ \mu $, $\lambda$, $\iota $, and $ \nu $. Therefore, applying the Aubin-Lions lemma yields that there exists $ (\vec u, h, A) $ satisfying \eqref{limit:regularity} and \eqref{limit:est-regularity},
such that, as $ (\mu,\lambda,\iota,\nu) \rightarrow (0^+,0^+,0^+,0^+) $,
\begin{equation}\label{limit:convergence}
    \begin{aligned}
    \vec u_{\mu,\lambda,\iota,\nu} & \rightarrow  \vec u & \text{in} \quad C(0,T^{**};H^3(\Omega)), \\
    h_{\mu,\lambda, \iota,\nu} & \rightarrow  h & \text{in} \quad C(0,T^{**};H^2(\Omega)), \\
    A_{\mu,\lambda, \iota,\nu} & \rightarrow  A & \text{in} \quad C(0,T^{**};H^2(\Omega)), \\
    (\vec u_{\mu,\lambda, \iota,\nu}, h_{\mu,\lambda, \iota,\nu},A_{\mu,\lambda, \iota,\nu}) &\buildrel\ast\over\rightharpoonup (\vec u, h, A) & \text{in} \quad L^\infty(0,T^{**};H^3(\Omega)), \\
    \vec u_{\mu,\lambda, \iota, \nu} & \rightharpoonup \vec u &\text{in}\quad L^2(0,T^{**};H^4(\Omega)), \\
    ( \dt \vec u_{\mu,\lambda, \iota,\nu}, \dt h_{\mu,\lambda, \iota,\nu}, \dt A_{\mu,\lambda, \iota,\nu} ) & \buildrel\ast\over\rightharpoonup (\dt \vec u, \dt h, \dt A) &\text{in} \quad L^\infty(0,T^{**};L^2(\Omega)),
    \end{aligned}
\end{equation}
and it is easy to verify that $ (\vec u, h, A) $ satisfies system \eqref{eq:sea-ice-app-2} in $ (0,T^{**}] $.

\subsection{Well-posedness of solutions for system \eqref{eq:sea-ice-app-2}}\label{subsec:well-posedness-app-2}

To deduce the well-posedness of solutions to system \eqref{eq:sea-ice-app-2},
it remains to establish the uniqueness and the continuous dependency of solutions on initial data. Indeed, this can be done following similar arguments as in Section \ref{subsec:contracting-mapping}, which we will sketch below.

Denote by $ (\vec u_\mrm j, h_\mrm j, A_\mrm j) $, $ \mrm j = 1,2 $, two solutions to system \eqref{eq:sea-ice-app-2} with initial data $ (\vec u_{\mrm{in,j}},h_{\mrm{in,j}}, A_{\mrm{in,j}}) $ within $ (0,T^{**}_{\mrm j}] $, $ \mrm j = 1, 2 $, as constructed above, respectively. In particular, \eqref{limit:regularity} and \eqref{limit:est-regularity} hold for $ (\vec u_\mrm j, h_\mrm j, A_\mrm j) $, $ \mrm j = 1, 2 $. Further, let $ \delta \vec u := \vec u_1 - \vec u_2,  \delta h:= h_1 - h_2, \delta A:= A_1 - A_2 $, and $  T_{12}^{**} := \min\lbrace T_1^{**}, T_2^{**} \rbrace \in (0,\infty) $. The triple $ (\delta \vec u, \delta h, \delta A) $ satisfies the following equations:
\begin{subequations}
	\begin{gather}
    \label{eq:u-sea-ice-diff}
    \begin{gathered}
        \rho_\mrm{ice} h_1 \dt \delta \vec u + \rho_\mrm{ice} \delta h \dt \vec u_2 = \dv \delta\mathbb S_{\varepsilon} - \nabla \delta p \\
        - \rho_\mrm{ice} h_1 \vec u_1 \cdot \nabla \delta \vec u - \rho_\mrm{ice} h_1 \delta \vec u \cdot \nabla \vec u_2 - \rho_\mrm{ice} \delta h \vec u_2 \cdot \nabla \vec u_2 + \delta \mathcal F,
    \end{gathered}\\
    \label{eq:h-sea-ice-diff}
    \dt \delta h + \dv(\delta h \vec u_1) + \dv ( h_2 \delta \vec u) = \delta \mathcal S_{h,\omega}, \\
    \label{eq:A-sea-ice-diff}
    \begin{gathered}
        \dt \delta A + \dv(\delta A \vec u_1) + \dv ( A_2 \delta\vec u) = \delta \mathcal S_{A,\omega} + \delta A \dv \vec u_1 \cdot \chi_{A_1}^\omega\\
        + A_2 \dv \delta \vec u \cdot \chi^\omega_{A_1} + A_2 \dv \vec u_2 \cdot \delta \chi_A^\omega.
    \end{gathered}
    \end{gather}
\end{subequations}

After taking the $ L^2 $-inner product of \eqref{eq:u-sea-ice-diff}, \eqref{eq:h-sea-ice-diff}, and \eqref{eq:A-sea-ice-diff} with $ 2 \delta \vec u $, $ 2 \delta h$ , and $ 2 \delta A $, respectively, and applying integration by parts in the resultant, one has
\begin{gather}
    \label{eq:diff-1}
    \begin{gathered}
        \dfrac{d}{dt} \norm{\rho_\mrm{ice}^{1/2} h_1^{1/2} \delta \vec u}{\Lnorm{2}}^2
        = \underbrace{- 2 \int \delta \mathbb S_\varepsilon: \nabla \delta \vec u\idx}_{\mathcal I_{14}}
        + \underbrace{\int \rho_\mrm{ice} \dt h_1 |\delta \vec u|^2 \idx}_{\mathcal I_{15}} \\
        \underbrace{- 2 \int \rho_\mrm{ice} \delta h \dt \vec u_2 \cdot \delta \vec u\idx}_{\mathcal I_{16}}
        + \underbrace{2 \int \delta p \, \dv \delta \vec u \idx}_{\mathcal I_{17}} + \underbrace{2 \int \delta\mathcal F \cdot \delta \vec u\idx}_{\mathcal I_{18}}
        \\ \underbrace{- 2 \int \rho_\mrm{ice} \bigl(h_1 \vec u_1 \cdot \nabla \delta \vec u + h_1 \delta \vec u\cdot \nabla \vec u_2 + \delta h \vec u_2 \cdot \nabla \vec u_2 \bigr) \cdot \delta \vec u\idx}_{\mathcal I_{19}}
    \end{gathered}\\
    \label{eq:diff-2}
    \begin{gathered}
        \dfrac{d}{dt} \norm{\delta h}{\Lnorm{2}}^2 = \underbrace{- \int \dv \vec u_1 |\delta h|^2 \idx}_{\mathcal I_{20}} \underbrace{- 2 \int \dv (h_2 \delta \vec u) \delta h\idx}_{\mathcal I_{21}} \\
        + \underbrace{2 \int \delta \mathcal S_{h,\omega} \delta h\idx}_{\mathcal I_{22}}, 
    \end{gathered}\\
    \label{eq:diff-3}
    \begin{gathered}
        \dfrac{d}{dt} \norm{\delta A}{\Lnorm{2}}^2 = \underbrace{- \int \dv \vec u_1 |\delta A|^2 \idx}_{\mathcal I_{23}} \underbrace{- 2 \int \dv (A_2 \delta \vec u) \delta A\idx}_{\mathcal I_{24}} \\
        + \underbrace{2 \int \delta \mathcal S_{A,\omega} \delta A\idx}_{\mathcal I_{25}} + \underbrace{2 \int \dv \vec u_1 \cdot \chi_{A_1}^\omega |\delta A|^2 \idx}_{\mathcal I_{26}} \\
        + \underbrace{2 \int A_2 \dv \delta \vec u \cdot \chi_{A_1}^\omega \delta A \idx}_{\mathcal I_{27}} + \underbrace{2 \int A_2 \dv \vec u_2 \cdot \delta \chi_{A}^\omega \delta A \idx}_{\mathcal I_{28}}.
    \end{gathered}
\end{gather}
Then it is straightforward to check that, thanks to the uniform bounds in \eqref{limit:est-regularity},
\begin{equation}\label{diffest-001}
    \begin{gathered}
        \sum_{15 \leq \mrm j \leq 28} \mathcal I_\mrm j \lesssim \norm{\delta\vec u,\delta h, \delta A}{\Lnorm{2}}^2 + \norm{\delta\vec u,\delta h, \delta A}{\Lnorm{2}}\norm{\nabla \delta \vec u}{\Lnorm{2}} .
    \end{gathered}
\end{equation}
To estimate $ \mathcal I_{14} $, we will have to investigate the monotonicity of $ \mathbb S_\varepsilon $, which is an important ingredient in our proof. Notice that
\begin{gather*}
    2 \delta \mathbb S_\varepsilon : \nabla \delta \vec u = \delta \biggl( p \dfrac{\nabla \vec u + \nabla \vec u^\top}{\sqrt{|\nabla \vec u + \nabla \vec u^\top|^2 + \varepsilon^2}}\biggr) : \delta (\nabla \vec u + \nabla \vec u^\top)  \\
    + 2 \delta \biggl( p \dfrac{\dv \vec u}{\sqrt{|\dv \vec u|^2 + \varepsilon^2}} \biggr) \delta \dv \vec u.
\end{gather*}
For $ \mrm D \vec u \in \lbrace \nabla \vec u + \nabla \vec u^\top , \dv \vec u \rbrace $, direct calculation yields that
\begin{equation}\label{id:monotone}
\begin{gathered}
    \delta \biggl( p \dfrac{\mrm D \vec u}{\sqrt{|\mrm D\vec u|^2 + \varepsilon^2}} \biggr) = \dfrac{1}{2} \biggl( \dfrac{p_1}{\sqrt{|\mrm D\vec u_2|^2 + \varepsilon^2}} + \dfrac{p_2}{\sqrt{|\mrm D\vec u_1|^2 + \varepsilon^2}} \biggr) \delta \mrm D \vec u  \\
    - \dfrac{1}{2} \dfrac{\bigl((\mrm D \vec u_1 + \mrm D \vec u_2) \cdot \delta \mrm D \vec u\bigr) \times ( p_1 \mrm D\vec u_1+ p_2 \mrm D\vec u_2)}{\sqrt{|\mrm D\vec u_1|^2 + \varepsilon^2}\sqrt{|\mrm D\vec u_2|^2 + \varepsilon^2}(\sqrt{|\mrm D\vec u_1|^2 + \varepsilon^2}+\sqrt{|\mrm D\vec u_2|^2 + \varepsilon^2})}\\
    + \dfrac{\delta p}{2}\biggl(\dfrac{\mrm D\vec u_1}{\sqrt{|\mrm D\vec u_1|^2 + \varepsilon^2}} +\dfrac{\mrm D\vec u_2}{\sqrt{|\mrm D\vec u_2|^2 + \varepsilon^2}} \biggr).
\end{gathered}
\end{equation}
Therefore
\begin{gather*}
    \delta \biggl( p \dfrac{\mrm D \vec u}{\sqrt{|\mrm D\vec u|^2 + \varepsilon^2}} \biggr) \cdot \delta \mrm D \vec u = \dfrac{\delta p}{2} \biggl(\dfrac{\mrm D\vec u_1}{\sqrt{|\mrm D\vec u_1|^2 + \varepsilon^2}} +\dfrac{\mrm D\vec u_2}{\sqrt{|\mrm D\vec u_2|^2 + \varepsilon^2}} \biggr) \cdot \delta \mrm D \vec u\\
    + \dfrac{1}{2} \dfrac{\mrm M}{\sqrt{|\mrm D\vec u_1|^2 + \varepsilon^2}\sqrt{|\mrm D\vec u_2|^2 + \varepsilon^2}(\sqrt{|\mrm D\vec u_1|^2 + \varepsilon^2}+\sqrt{|\mrm D\vec u_2|^2 + \varepsilon^2})},
\end{gather*}
with
\begin{gather*}
    \mrm M := \bigl(p_1\sqrt{|\mrm D u_1|^2 + \varepsilon^2} + p_2 \sqrt{|\mrm D u_2|^2 + \varepsilon^2} \bigr) \bigl( \sqrt{|\mrm D \vec u_1|^2 + \varepsilon^2} + \sqrt{|\mrm D \vec u_2|^2 + \varepsilon^2} \bigr) \\
     \times | \delta \mrm D \vec u|^2
    - \bigl((\mrm D \vec u_1 + \mrm D \vec u_2) \cdot \delta \mrm D \vec u\bigr) \times \bigl(( p_1 \mrm D\vec u_1+ p_2 \mrm D\vec u_2) \cdot \delta \mrm D \vec u\bigr) \\
    \geq C_{\underline h,\mathfrak C_\mrm{in}}\varepsilon | \delta \mrm D \vec u|^2,
\end{gather*}
for some constant $ C_{\underline h,\mathfrak C_\mrm{in}} \in (0,\infty) $ depending on $ \underline h $ and $ \mathfrak C_\mrm{in}$.
Therefore, one can derive that
\begin{equation}\label{diffest-002}
\begin{gathered}
    \mathcal I_{14} \lesssim - C_{\varepsilon, \underline h,\mathfrak C_\mrm{in}} \bigl( \norm{\nabla \delta \vec u + \nabla \delta \vec u^\top}{\Lnorm{2}}^2 + \norm{\dv \delta \vec u}{\Lnorm{2}}^2 \bigr) \\
    + \norm{\delta h, \delta A}{\Lnorm{2}}\norm{\nabla \delta \vec u}{\Lnorm{2}},
\end{gathered}
\end{equation}
for some constant $ C_{\varepsilon,\underline h, \mathfrak C_\mrm{in}} \in (0,\infty) $ depending on $ \varepsilon $, $ \underline h $, and $ \mathfrak C_\mrm{in}$.
In addition, using integration by parts, one can derive that,
\begin{equation}\label{diffest-003}
\norm{\nabla \delta \vec u}{\Lnorm{2}}^2 \lesssim \norm{\nabla \delta \vec u + \nabla \delta \vec u^\top}{\Lnorm{2}}^2 + \norm{\dv \delta \vec u}{\Lnorm{2}}^2.
\end{equation}

Consequently,
after substituting \eqref{diffest-001}, \eqref{diffest-002}, and \eqref{diffest-003} into \eqref{eq:diff-1}, \eqref{eq:diff-2}, and \eqref{eq:diff-3}, summing up the results, and applying Young's inequality, one can conclude that
\begin{equation*}
    \dfrac{d}{dt} \norm{\rho_\mrm{ice}^{1/2} h_1^{1/2} \delta \vec u, \delta h, \delta A}{\Lnorm{2}}^2 \leq C_{\mathfrak C_{\mrm in}} \norm{\rho_\mrm{ice}^{1/2} h_1^{1/2} \delta \vec u, \delta h, \delta A}{\Lnorm{2}}^2,
\end{equation*}
which, after applying Gr\"onwall's inequality, yields
\begin{equation}\label{diffest-101}
    \sup_{0\leq s\leq T_{12}^{**}}\norm{\delta \vec u(s), \delta h(s), \delta A(s)}{\Lnorm{2}}^2 \leq C_{\mathfrak C_\mrm{in}} \norm{\delta u_\mrm{in}, \delta h_\mrm{in}, \delta A_\mrm{in}}{\Lnorm{2}}^2,
\end{equation}
with some constant $ C_{\mathfrak C_\mrm{in}} \in(0,\infty) $, depending on the initial data. The uniqueness and the continuous dependence on initial data of solutions to system \eqref{eq:sea-ice-app-2} follow from \eqref{diffest-101}.

\section*{Acknowledgement}
XL and MT gratefully acknowledge the partial funding by the Deutsche
Forschungsgemeinschaft (DFG) through project AA2-9 \emph{Variational methods for viscoelastic flows and gelation} within MATH$^{+}$. MT also gratefully acknowledges the partial funding by the DFG through project
C09 \emph{Dynamics of rock dehydration on multiple scales} (project number 235221301) within CRC 1114 \emph{Scaling Cascades in Complex Systems}. Moreover  XL and EST are thankful for the kind hospitality of Freie Universit\"at Berlin where part of this work was done and partially supported by the Einstein Stiftung/Foundation - Berlin, through the Einstein Visiting Fellow Program. 
	EST and XL would also like to thank the Isaac Newton Institute for Mathematical Sciences for support and hospitality during the programme TUR when part of this work was undertaken. This work was supported by EPSRC Grant Number EP/R014604/1.
	XL's work was partially supported by a grant from the Simons Foundation, during his visit to the Isaac Newton Institute. The authors also thank the reviewers for the helpful comments during the submission of this work. 

\bibliographystyle{plain}

\end{document}